# On Subgroups of Tetrahedron Groups

Ma. Louise N. De Las Peñas[1], Rene P. Felix[2] and Glenn R. Laigo[3]

[1]*Mathematics Department, Ateneo de Manila University, Loyola Heights, Quezon City, Philippines. mlp@math.admu.edu.ph*

[2]*Institute of Mathematics, University of the Philippines – Diliman, Diliman, Quezon City, Philippines. rene@upd.edu.ph*

[3]*Mathematics Department, Ateneo de Manila University, Loyola Heights, Quezon City, Philippines. glaigo@ateneo.edu*

**Abstract**. We present a framework to determine subgroups of tetrahedron groups and tetrahedron Kleinian groups, based on tools in color symmetry theory.

## 1 Introduction

In [6, 7] the determination of the subgroups of two-dimensional symmetry groups was facilitated using a geometric approach and applying concepts in color symmetry theory. In particular, these works derived low-indexed subgroups of a triangle group **abc*, a group generated by reflections along the sides of a given triangle $\Delta$ with interior angles $\pi/a$, $\pi/b$ and $\pi/c$ where $a$, $b$, $c$ are integers $\geq 2$. The elements of **abc* are symmetries of the tiling formed by copies of $\Delta$. The subgroups of **abc* were calculated from colorings of this given tiling, assuming a corresponding group of color permutations.

In this work, we pursue the three dimensional case and extend the problem to determine subgroups of tetrahedron groups, groups generated by reflections along the faces of a Coxeter tetrahedron, which are groups of symmetries of tilings in space by tetrahedra. The study will also address determining subgroups of other types of three dimensional symmetry groups such as the tetrahedron Kleinian groups.

## 2 Coxeter tetrahedra, tetrahedron groups and tetrahedron Kleinian groups

Let $\mathbf{X^3} = \mathbf{S}^3$, $\mathbf{E}^3$, or $\mathbf{H}^3$ denote either the spherical, Euclidean or hyperbolic 3-space. Consider a tetrahedron $\mathtt{t}$ in $\mathbf{X^3}$ all of whose dihedral angles are submultiples of $\pi$. $\mathtt{t}$ is called a *Coxeter tetrahedron* and if its dihedral angles are given by $\pi/p$, $\pi/q$, $\pi/r$, $\pi/s$, $\pi/t$, $\pi/u$ we denote $\mathtt{t}$ by $[p, q, r, s, t, u]$. Five of the known Coxeter tetrahedra are spherical: [3, 3, 3, 2, 2, 2], [4, 3, 3, 2, 2, 2], [3, 4, 3, 2, 2, 2], [5, 3, 3, 2, 2, 2], and [3, 3, 2, 3, 2, 2]; and three are Euclidean: [4, 3, 4, 2, 2, 2], [3, 4, 2, 3, 2, 2], and [3, 3, 3, 2, 2, 3][2]. In the hyperbolic case, there are 9 compact and 23 non-compact Coxeter tetrahedra of finite volume [5, 13]. The compact types have no ideal vertices or vertices at infinity; while the 23 non-compact ones have at least one ideal vertex.

Given a Coxeter tetrahedron $\mathtt{t}$ in $\mathbf{X^3}$, the reflections in its faces generate a group called a *tetrahedron group* or a *Coxeter tetrahedron group*. The subgroup of index 2 in the tetrahedron group consisting of orientation preserving isometries is referred to as a *tetrahedron Kleinian group*. We address in this paper the problem of determining the subgroups of tetrahedron groups and tetrahedron

Kleinian groups corresponding to hyperbolic Coxeter tetrahedra of finite volume. Studies have been carried out on the subgroup structure of tetrahedron groups and tetrahedron Kleinian groups associated with spherical and Euclidean Coxeter tetrahedra, but a lot is yet to be discovered for the hyperbolic case.

In our work, the determination of the subgroups of tetrahedron groups and tetrahedron Kleinian groups is facilitated by the link of group theory and color symmetry theory, as will be discussed in the next section.

## 3 Colorings of three-dimensional tilings by tetrahedra

We begin by first explaining the setting we will be working with in studying tetrahedron groups and their subgroups, and the notation we will use in the results that follow. Consider a Coxeter tetrahedron $\mathtt{t} := [p, q, r, s, t, u]$ of finite volume with vertices $P_0$, $P_1$, $P_2$, $P_3$ and dihedral angles $\pi/p$, $\pi/q$, $\pi/r$, $\pi/s$, $\pi/t$, $\pi/u$, as shown in Fig. 1, lying on either $\mathbf{X^3} = \mathbf{S}^3$, $\mathbf{E}^3$ or $\mathbf{H}^3$. Repeatedly reflecting $\mathtt{t}$ in its faces results in a tiling $\mathscr{T}$ of $\mathbf{X^3}$ by copies of $\mathtt{t}$. Let $P$, $Q$, $R$, $S$ denote reflections along the respective faces $P_1P_2P_3$, $P_0P_2P_3$, $P_0P_1P_3$, $P_0P_1P_2$ of $\mathtt{t}$. The group $H$ generated by $P$, $Q$, $R$, $S$ is called a tetrahedron group. The subgroup of the tetrahedron group consisting of orientation preserving isometries, and in particular generated by $PQ$, $QR$, and $RS$ is the tetrahedron Kleinian group $K$.

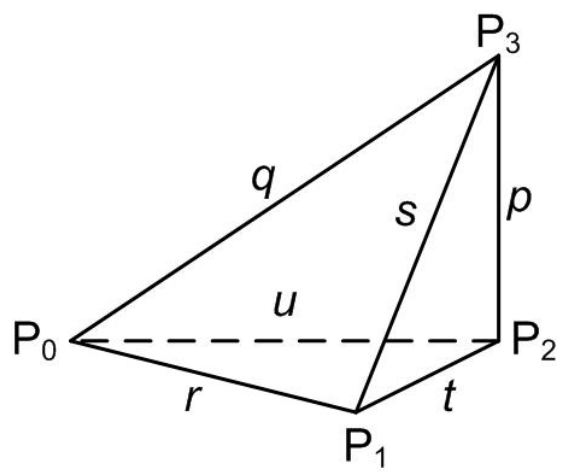

Fig.1. A tetrahedron $\mathtt{t}$ with vertices $P_0$, $P_1$, $P_2$, $P_3$. A label $k$ at an edge of $\mathtt{t}$ means that the dihedral angle of $\mathtt{t}$ at this edge is $\pi/k$.

As in [26], we have the following result on the tetrahedron group $H$:

**Lemma 1.** *The tiling $\mathscr{T}$ is the H-orbit of $\mathtt{t}$, that is, $\mathscr{T} = \{h\mathtt{t} : h \in H\}$. Moreover, $\mathtt{t}$ forms a fundamental region for H and $Stab_H(\mathtt{t}) = \{h \in H : h\mathtt{t} = \mathtt{t}\}$ is the trivial group $\{e\}$.*

By this lemma we mean that if $h \in H$ and $h\mathtt{t}$ is the image of $\mathtt{t}$ under $h$ then the union of these images as $h$ varies over the elements of $H$ is $\mathbf{X}^3$ and, moreover, if $h_1$, $h_2 \in H$, $h_1 \neq h_2$, the respective interiors of $h_1\mathtt{t}$ and $h_2\mathtt{t}$ are disjoint. Each tetrahedron in $\mathscr{T}$ is the image of $\mathtt{t}$ under a uniquely determined $h \in H$.

Suppose $M$ is a subgroup of $H$ and $O = M\mathtt{t} = \{m\mathtt{t} : m \in M\}$, the $M$-orbit of $\mathtt{t}$. Then $Stab_M(\mathtt{t}) = \{e\}$ and $M$ acts transitively on $O$. Consequently, there is a one-to-one correspondence between $M$ and $O$ given by $m \to m\mathtt{t}$, $m \in M$. The action of $M$ on $O$ is regular, where $m' \in M$ acts on $m\mathtt{t} \in O$ by sending it to its image under $m'$.

In this paper, the study of the subgroups of the tetrahedron group $H$ and tetrahedron Kleinian group $K$ is addressed via tools in color symmetry theory, and in particular, using colorings of $\mathscr{T}$ and $K\mathtt{t}$, respectively.

We give the definition of a coloring of the $M$-orbit of $\mathtt{t}$ for a subgroup $M$ of $H$ as follows:

Consider $O = M\mathtt{t} = \{m\mathtt{t} : m \in M\}$, $M$ a subgroup of $H$. If $C = \{c_1, c_2, \ldots, c_n\}$ is a set of $n$ colors, an onto function $f : O \rightarrow C$ is called an *n-coloring* of $O$. To each $m\mathtt{t} \in O$ is assigned a color in $C$. The coloring determines a partition $\mathcal{P} = \{f^{-1}(c_i) : c_i \in C\}$ where $f^{-1}(c_i)$ is the set of elements of $O$ assigned color $c_i$. Equivalently, we may think of the coloring as a partition of $O$.

We look at $n$-colorings of $O$ for which the group $M$ has a transitive action on the set $C$ of $n$ colors. We will refer to such a coloring as an $M$-transitive $n$-coloring of $O$. To arrive at such a coloring, consider the following lemma.

**Lemma 2.** *Let $\mathcal{P} = \{X_1, X_2, \ldots, X_n\}$ be a partition of $O$ into n sets. Suppose M acts transitively on $\mathcal{P}$, then $\mathcal{P} = \{mL\mathtt{t} : m \in M\}$ for some subgroup L of M of index n.*

***Proof***

Assume $\mathtt{t} \in X_1$ and let $L = Stab_M(X_1)$. Using the fact that $M$ acts transitively on $O$ and $M$ acts transitively on $\mathcal{P}$, it follows that $X_1 = L\mathtt{t}$ and each set in the partition $\mathcal{P}$ is equal to $mX_1$ for some $m \in M$. Thus, $\mathcal{P} = \{mL\mathtt{t} : m \in M\}$. Moreover, by the orbit-stabilizer theorem, $n = |\mathcal{P}| = [M : Stab_M(X_1)] = [M : L]$.■

Let $L$ be a subgroup of $M$ of index $n$ and $\{m_1, m_2, \ldots, m_n\}$ be a complete set of left coset representatives of $L$ in $M$. Define the coloring $f : O = M\mathtt{t} \rightarrow C = \{c_1, c_2, \ldots, c_n\}$ by $f(x) = c_i$ if $x \in m_iL\mathtt{t}$. Then $f^{-1}(c_i) = m_iL\mathtt{t}$ and $\mathcal{P} = \{f^{-1}(c_1), f^{-1}(c_2), \ldots, f^{-1}(c_n)\} = \{m_1L\mathtt{t}, m_2L\mathtt{t}, \ldots, m_nL\mathtt{t}\}$. The group $M$ acts transitively on $\mathcal{P}$ with $m \in M$ sending $m_iL\mathtt{t}$ to its image $mm_iL\mathtt{t}$. Consequently, we get a transitive action of $M$ on $C$ by defining for $m \in M$, $mc_i = c_j$ if and only if $mf^{-1}(c_i) = f^{-1}(c_j)$. The $n$-coloring of $O$ determined by $f$ is $M$- transitive.

Our main result which is a consequence of the preceding discussions and forms the basis of our calculations is given as follows.

**Theorem.** *Let M be a subgroup of* $H = \langle P, Q, R, S\rangle$ and $O = \{m\mathtt{t} : m \in M\}$ *the M-orbit of* $\mathtt{t}$.

(i) *Suppose L is a subgroup of M of index n. Let $\{m_1, m_2, \ldots, m_n\}$ be a complete set of left coset representatives of L in M and $\{c_1, c_2, \ldots, c_n\}$ a set of n colors. The assignment $m_iL\mathtt{t} \rightarrow c_i$ defines an n-coloring of O which is M-transitive.*

(ii) *In an M-transitive n-coloring of O, the elements of M which fix a specific color in the colored set O form a subgroup of M of index n.*

**Remark:** Note that given a subgroup $L$ of $M$ of index $n$ and a set of $n$ colors $\{c_1, c_2, \ldots, c_n\}$, then there correspond $(n-1)!$ $M$-transitive $n$-colorings of $O$ with $L$ fixing $c_1$. In an $M$-transitive $n$-coloring of $O$ with $L$ fixing $c_1$, the set of tetrahedra $L\mathtt{t}$ is assigned $c_1$ and the remaining $n-1$ colors are distributed among the $m_iL\mathtt{t}$, $m_i \notin L$.

## 4 Index 2, 3 and 4 subgroups of tetrahedron groups

In applying the given theorem to derive the index 2, 3 and 4 subgroups of the tetrahedron group $H$, we construct respectively, two, three and four colorings of $\mathcal{T}$ where all elements of $H$ effect color permutations and $H$ acts transitively on the set of colors. For such an $n$-coloring of $\mathcal{T}$, $n \in \{2, 3, 4\}$, a homomorphism $\pi: H \to S_n$ is defined. The group $H$ is generated by $P$, $Q$, $R$ and $S$, thus $\pi$ is completely determined when $\pi(P)$, $\pi(Q)$, $\pi(R)$ and $\pi(S)$ are specified. To construct a coloring, we consider a fundamental region for $H$ and assign to it color $c_1$. The arguments to obtain 2-, 3- or 4-colorings of $\mathcal{T}$ where each of the generators $P$, $Q$, $R$, $S$ permutes the colors, that is, the associated color permutation to a generator is of order 1 or 2, are as follows:

(*i*) In constructing 2-colorings of $\mathcal{T}$, each of $P$, $Q$, $R$, $S$ may be assigned either the identity permutation (1) or the permutation (12). The possible 2-colorings of $\mathcal{T}$ are listed in Table 1.
(*ii*) In arriving at 3-colorings of $\mathcal{T}$, this time each of $P$, $Q$, $R$, $S$ may be assigned either the identity permutation or either of the permutations (12), (13) or (23). Now, since we assume $H$ acts transitively on the set of 3 colors, we consider the assignments that bring forth a $\pi(H)$ which is a transitive subgroup of $S_3$. Note that there are two transitive subgroups of $S_3$ namely $S_3$ and the cyclic group $Z_3$. It is not possible for $\pi(H)$ to be $Z_3$ since this group cannot be generated by elements of order two. Thus we only consider the case when $\pi(H)$ is $S_3$. In order to generate $S_3$, two of the three permutations (12), (13) or (23) must be associated with two from among the generators $P$, $Q$, $R$ and $S$. The remaining third and fourth generators are associated either with the identity, with one of the two already assigned earlier or with the third permutation of order 2. The possible 3-colorings that will yield distinct subgroups of index 3 distinct up to conjugacy in $H$ are given in Table 2.
(*iii*) In constructing 4-colorings of $\mathcal{T}$ that will yield index 4 subgroups of $H$, each of $P$, $Q$, $R$, $S$ is assigned either the identity permutation; a 2-cycle or a product of two disjoint 2-cycles: (12), (13), (14), (23), (24), (34), (12)(34),(13)(24) or (14)(23). We only look at the situations that arise when $\pi(H)$ is a transitive subgroup of $S_4$. In particular, we consider the cases when $\pi(H)$ is either $S_4$, $D_4$ or $V = \{e, (12)(34), (13)(24), (14)(23)\}$, a Klein 4 group. The possible $H$-transitive four-colorings of $\mathcal{T}$ that can be constructed when $\pi(H) \cong V$, $D_4$ and $S_4$ are given respectively in Tables 3, 4 and 5. These colorings will give rise to distinct subgroups of index 4 distinct up to conjugacy in $H$.

To illustrate the process in determining subgroups of tetrahedron groups, let us consider the tetrahedron group $H$ generated by the reflections $P$, $Q$, $R$, $S$ along the faces of the Coxeter tetrahedron $\mathtt{t}_{10}$ := [3, 3, 6, 2, 2, 2] (refer to Table 7 for the complete list of hyperbolic Coxeter tetrahedra). $\mathtt{t}_{10}$ is a tetrahedron with one vertex at infinity and four right angled triangles for its faces. The symmetry group *G* of the tiling $\mathscr{T}$ by copies of $\mathtt{t}_{10}$ is $H$ = < *P, Q, R, S* > with defining relations $P^2 = Q^2 = R^2 = S^2 = (PQ)^3 = (QR)^3 = (RS)^6 = (PR)^2 = (PS)^2 = (QS)^2 = e$.

There are only three index 2 subgroups of *H*, namely, <*P*, *Q*, *R*, *SRS*>, <*QP, RP, S*> and <*PQ, QR, RS*> corresponding respectively to the 2-colorings given in Table 1 line nos 4, 14 and 15. On the other hand *H* has only one index 3 subgroup distinct up to conjugacy in *H*, the group <*S, QPQ, QRQ, RP*> corresponding to the 3-coloring listed in Table 2 line no. 14. There are two index 4 subgroups of *H*: <*QP*, *RP*, *SRSP*> and <*R, S, QPQ, QRSRQ*> corresponding respectively to the 4-colorings given in Table 3 line no 32 and Table 5 line no 92.

Note that the index 2 subgroup listed in Table 1 no. 4 with six generators may be actually generated by four generators, namely, *P*, *Q*, *R*, *SRS*; the index 2 subgroup given in Table 1 no. 14 with four generators may be generated by three, given by *QP, RP, S* and the index 3 subgroup of *H* listed in Table 2, line no 14 with 7 generators, may be generated by four generators: *S, QPQ, QRQ, RP*. Also, the index 4 subgroup may have three generators: *QP*, *RP*, *SRSP* instead of five given in Table 3 line no 32 and the index 4 subgroup may have 4 generators: *R, S, QPQ, QRSRQ* instead of 10 as listed in Table 5 line no 92. The involutions *QS*, *PS* and *PR* give rise respectively to the relations $QSQ = S$, $PSP = S$ and $PRP = R$; the order 3 rotations *PQ* and *QR* result in $PQP = QPQ$ and $RQR = QRQ$.

Note that we do not consider the permutation assignments to *P*, *Q*, *R*, *S* that will result in permutations corresponding to *PQ* and *QR* that are 2- or 4-cycles; *PR*, *PS*, and *QS* that are 3-cycles and *RS* that is a 4-cycle. This is due to the fact that *PQ* and *QR* are order 3 elements; *PR*, *PS*, and *QS* are order 2 elements and *RS* has order 6.

The tetrahedron group *H* corresponding to $\mathtt{t}_{10}$ has an interesting subgroup structure. Some of its subgroups appear in existing literature [12, 25]. For instance, its index 2 subgroup, the tetrahedron Kleinian group <*PQ*, *QR*, *RS*>, is the extended Eisenstein modular group $P\overline{S}L_2(E)$.The index 4 subgroup of *H* generated by the three rotations *QP*, *RP*, and *SRSP* form a subgroup of the special linear group $SL_2(E)$ and is the Eisenstein modular group $PSL_2(E)$. The second index 4 subgroup generated by the reflections *R*, *S*, *QPQ*, and *QRSRQ* is the symmetry group of the regular honeycomb {3, 6, 3} and has fundamental region a Coxeter tetrahedron of type $\mathtt{t}_{15}$.

## 5 Index 2, 3 and 4 subgroups of tetrahedron Kleinian groups

In determining the subgroups of tetrahedron Kleinian groups, we carry over the same setting given in the previous section, that is, we consider a Coxeter tetrahedron $\mathtt{t}$ := [*p, q, r, s, t*] of finite volume and the three-dimensional tiling $\mathscr{T}$ which results by reflecting $\mathtt{t}$ in its faces. The group generated by the reflections *P, Q, R, S* in the faces of $\mathtt{t}$, is the tetrahedron group *H*. The subgroup of the tetrahedron group generated by *PQ*, *QR*, and *RS* is the tetrahedron Kleinian group *K*.

In applying our theorem to derive the subgroups of *K,* we consider the *K*-orbit of $\mathtt{t}$, namely $O = K\mathtt{t} = \{k\mathtt{t} : k \in K\}$ and construct *K*-transitive colorings of $O$. In determining respectively, the index 2, 3 and 4 subgroups of *K*, we construct two, three and four colorings of $O$. For an $n$-coloring of $O$, $n \in \{2, 3, 4\}$, a homomorphism $\pi': K \to S_n$ is defined. The group $K$ is generated by $PQ$, $QR$ and $RS$, thus $\pi'$ is completely determined when $\pi'(PQ)$, $\pi'(QR)$, $\pi'(RS)$ are specified. To obtain 2-, 3- or 4-colorings of $O$ where each of the generators $PQ$, $QR$, $RS$ permutes the colors, the orders of the color permutations associated to $PQ$, $QR$, $RS$ should be divisors respectively of the orders of $PQ$, $QR$, $RS$. Construction of tables similar to Tables 1-5, this time pertaining to assignments of permutations to $PQ$, $QR$ and $RS$, will yield transitive two, three and four colorings of $O$.

To illustrate the above ideas, let us consider the tetrahedron Kleinian subgroup $K = \langle PQ, QR, RS \rangle$ corresponding to the Coxeter tetrahedron $\mathtt{t}_{10}$ with defining relations $(PQ)^3 = (QR)^3 = (RS)^6 = (PR)^2 = (QS)^2 = (PS)^2 = e$. The colorings that will yield low index subgroups of $K$ are listed in Table 6.

To obtain 2-colorings of $O$ that will yield index 2 subgroups of $K$, we assign to $PQ$, $QR$, $RS$ either the identity permutation (1) or the 2-cycle (12). Observe that $PQ$ and $QR$ have odd orders so that $PQ$ and $QR$ have to be assigned the identity permutation. There is only one such permutation assignment that satisfies these conditions, giving rise to one index 2 subgroup.

In obtaining 3-colorings of $O$ that will give rise to index 3 subgroups of $K$, we assign to $PQ$, $QR$, $RS$ color permutations that will generate a transitive subgroup of $S_3$. We will only consider the permutation assignments that bring forth a $\pi'(K)$ isomorphic to $Z_3$ or $S_3$. Moreover, $PQ$ and $QR$ are order 3 elements so that $PQ$ and $QR$ can only be assigned either the identity or a 3-cycle. Also we only consider the permutation assignments to $PQ$, $QR$, $RS$ resulting in permutations corresponding to $PR$, $PS$, and $QS$ that are either the identity or a 2-cycle, since $PR$, $PS$ and $QS$ are involutions. With these conditions, we arrive at one index 3 subgroup of $K$ up to conjugacy in $K$, where $\pi'(K)$ is $S_3$. There are no index 3 subgroups of $K$ where $\pi'(K)$ is $Z_3$ for this example.

In determining the 4-colorings of $O$ that will correspond to index 4 subgroups of $K$, we assign to $PQ$, $QR$, $RS$ color permutations that will yield any transitive subgroup of $S_4$. We will consider the permutation assignments that bring forth a $\pi'(K)$ isomorphic to $V = \{e, (12)(34), (13)(24), (14)(23)\}$, $Z_4$, $D_4$, $A_4$, or $S_4$. As explained above, $PQ$ and $QR$ can only be assigned either the identity or a 3-cycle; $RS$ has order 6 and cannot be assigned a 4-cycle. Moreover, we consider the permutation assignments to $PQ$, $QR$, $RS$ resulting in permutations corresponding to $PR$, $PS$, and $QS$ that are either the identity or a 2-cycle, or a product of disjoint 2-cycles. With these conditions, we arrive at one index 4 subgroup of $K$, up to conjugacy in $K$, where $\pi'(K)$ is $S_4$. There are no index 4 subgroups of $K$ where $\pi'(K)$ is Klein 4, $Z_4$, $D_4$ or $A_4$ for this example.

## 6 Low index subgroups of tetrahedron groups and tetrahedron Kleinian groups

In this part of the paper, we present the number of low index subgroups of the tetrahedron groups and tetrahedron Kleinian groups arising from the Coxeter tetrahedra of finite volume. The results for the number of index 2, 3, 4 subgroups of these groups, up to conjugacy, are shown in Table 7. In column 1, we list the compact tetrahedra ($\mathtt{t}_1 - \mathtt{t}_9$) and the non-compact ones ($\mathtt{t}_{10} - \mathtt{t}_{32}$). The number of vertices at infinity is also given, for each non-compact tetrahedron.

## 7 Future Outlook

In this study, a method is provided which facilitates the enumeration of the subgroups of a group of symmetries of a three dimensional tiling by a Coxeter tetrahedron. The approach we have presented here allows one to study three dimensional Coxeter groups and their subgroups from a geometric perspective using concepts in color symmetry theory. In our work, we focused on deriving the index 2, 3 and 4 subgroups of hyperbolic tetrahedron groups and tetrahedron Kleinian groups realized by Coxeter tetrahedra of finite volume. In coming up with the list of subgroups, most of the information pertaining to the subgroups relied on the construction of a table of permutation assignments.

As a continuation of this study, it would be interesting to look at the subgroup structure of the other subgroups of a tetrahedron group using the approach presented here. The Picard group for example, which is an index 4 subgroup of the tetrahedron group associated with $\mathtt{t}_{11}$ [17] has a very interesting subgroup structure. The Picard group has especially been of great interest in many fields of mathematics, for example, number theory and automorphic function theory. Particular examples of subgroups of this group have been studied in existing literature [9, 12, 21, 24, 25]. Continuing the study of the subgroup structure of the tetrahedron group will also include determining its two dimensional subgroups. For example, the Picard group contains the modular group $PSL(2,Z)$ as a subgroup [8]. It would also be interesting to determine the classes of triangle fuchsian groups contained in the tetrahedron Kleinian groups.

The study of tetrahedra with infinite volume is still at its infancy, and not much has been said about its corresponding tetrahedron groups [1, 14]. Studying these groups and their subgroups would be worth looking at. There are initial studies on truncated tetrahedra, which are polyhedra with some or all of its vertices cut off arising from tetrahedra of infinite volume. When there are $m$ truncations, a corresponding group arises generated by $m + 4$ reflections [10, 18]. It would also be interesting to look at the subgroup structure of these types of reflection groups.

Finally, another avenue in the study of hyperbolic Coxeter groups and their subgroup structure would be to extend the approach given in this work to look at the higher dimensional cases. For example, there are 5 compact and 9 non-compact hyperbolic simplices in 4 dimensions; and 12 non-compact hyperbolic simplices in 5 dimensions. There are also existing non-compact hyperbolic simplices of finite volume in dimensions 6 to 9 [13]. The next step would be to determine the subgroups of the Coxeter groups corresponding to these simplices.

*Acknowledgements.* Louise De Las Peñas acknowledges the support from the FEBTC- David G. Choa Professorial Chair Faculty Grant, Ateneo de Manila University. Glenn R. Laigo is grateful for a grant from the Commission on Higher Education – Higher Education Development Program. It is a pleasure to acknowledge Marjorie Senechal for suggestions on the problem.

| no | *P* | *Q* | *R* | *S* | *PQ* $\pi/p$ | *QR* $\pi/q$ | *RS* $\pi/r$ | *PR* $\pi/s$ | *PS* $\pi/t$ | *QS* $\pi/u$ | generators for the subgroup fixing color 1 |
|---|---|---|---|---|---|---|---|---|---|---|---|
| 1 | (12) | (1) | (1) | (1) | (12) | (1) | (1) | (12) | (12) | (1) | *Q, R, S, PQP, PRP, PSP* |
| 2 | (1) | (12) | (1) | (1) | (12) | (12) | (1) | (1) | (1) | (12) | *P, R, S, QPQ, QRQ, QSQ* |
| 3 | (1) | (1) | (12) | (1) | (1) | (12) | (12) | (12) | (1) | (1) | *P, Q, S, RPR, RQR, RSR* |
| 4 | (1) | (1) | (1) | (12) | (1) | (1) | (12) | (1) | (12) | (12) | *P, Q, R, SPS, SQS, SRS* |
| 5 | (1) | (1) | (12) | (12) | (1) | (12) | (1) | (12) | (12) | (12) | *P, Q, SR, RPR, RQR* |
| 6 | (1) | (12) | (1) | (12) | (12) | (12) | (12) | (1) | (12) | (1) | *P, R, SQ, QPQ, QRQ* |
| 7 | (1) | (12) | (12) | (1) | (12) | (1) | (12) | (12) | (1) | (12) | *P, RQ, S, QPQ, QSQ* |
| 8 | (12) | (12) | (1) | (1) | (1) | (12) | (1) | (12) | (12) | (12) | *QP, R, S, PRP, PSP* |
| 9 | (12) | (1) | (12) | (1) | (12) | (12) | (12) | (1) | (12) | (1) | *Q, RP, S, PQP, PSP* |
| 10 | (12) | (1) | (1) | (12) | (12) | (1) | (12) | (12) | (1) | (12) | *Q, R, SP, PQP, PRP* |
| 11 | (1) | (12) | (12) | (12) | (12) | (1) | (1) | (12) | (12) | (1) | *P, RQ, SQ, QPQ* |
| 12 | (12) | (1) | (12) | (12) | (12) | (12) | (1) | (1) | (1) | (12) | *Q, RP, SP, PQP* |
| 13 | (12) | (12) | (1) | (12) | (1) | (12) | (12) | (12) | (1) | (1) | *QP, R, SP, PRP* |
| 14 | (12) | (12) | (12) | (1) | (1) | (1) | (12) | (1) | (12) | (12) | *QP, RP, S, PSP* |
| 15 | (12) | (12) | (12) | (12) | (1) | (1) | (1) | (1) | (1) | (1) | *PQ, QR, RS* |

Table 1. *H*-transitive colorings that will give rise to index 2 subgroups of *H*.

| no | *P* | *Q* | *R* | *S* | *PQ* $\pi/p$ | *QR* $\pi/q$ | *RS* $\pi/r$ | *PR* $\pi/s$ | *PS* $\pi/t$ | *QS* $\pi/u$ | generators for the subgroup fixing color 1 |
|---|---|---|---|---|---|---|---|---|---|---|---|
| 1 | (12) | (13) | (1) | (1) | (132) | (13) | (1) | (12) | (12) | (13) | *R, S, QPQ, QRQ, QSQ, PQP, PRP, PSP* |
| 2 | (12) | (1) | (13) | (1) | (12) | (13) | (13) | (132) | (12) | (1) | *Q, S, RPR, RQR, RSR, PQP, PRP, PSP* |
| 3 | (12) | (1) | (1) | (13) | (12) | (1) | (13) | (12) | (132) | (13) | *Q, R, SPS, SQS, SRS, PQP, PRP, PSP* |
| 4 | (1) | (12) | (13) | (1) | (12) | (132) | (13) | (13) | (1) | (12) | *P, S, RPR, RQR, RSR, QPQ, QRQ, QSQ* |
| 5 | (1) | (12) | (1) | (13) | (12) | (12) | (13) | (1) | (13) | (132) | *P, R, SPS, SQS, SRS, QPQ, QRQ, QSQ* |
| 6 | (1) | (1) | (12) | (13) | (1) | (12) | (132) | (12) | (13) | (13) | *P, Q, SPS, SQS, SRS, RPR, RQR, RSR* |
| 7 | (12) | (12) | (1) | (13) | (1) | (12) | (13) | (12) | (132) | (132) | *R, QP, SPS, SQS, SRS, PRP, PSP* |
| 8 | (12) | (1) | (12) | (13) | (12) | (12) | (132) | (1) | (132) | (13) | *Q, RP, SPS, SQS, SRS, PQP, PSP* |
| 9 | (12) | (1) | (13) | (12) | (12) | (13) | (123) | (132) | (1) | (12) | *Q, SP, RPR, RQR, RSR, PQP, PRP* |
| 10 | (1) | (12) | (12) | (13) | (12) | (1) | (132) | (12) | (13) | (132) | *P, RQ, SPS, SQS, SRS, QPQ, QSQ* |
| 11 | (1) | (12) | (13) | (12) | (12) | (132) | (123) | (13) | (12) | (1) | *P, SQ, RPR, RQR, RSR, QPQ, QRQ* |
| 12 | (1) | (13) | (12) | (12) | (13) | (12) | (1) | (12) | (12) | (123) | *P, RS, RPR, RQR, QPQ, QRQ, QSQ* |
| 13 | (12) | (12) | (13) | (1) | (1) | (132) | (13) | (132) | (12) | (12) | *S, QP, RPR, RQR, RSR, PRP, PSP* |
| 14 | (12) | (13) | (12) | (1) | (132) | (123) | (12) | (1) | (12) | (13) | *S, RP, QPQ, QRQ, QSQ, PQP, PSP* |
| 15 | (12) | (13) | (1) | (12) | (132) | (13) | (12) | (12) | (1) | (123) | *R, SP, QPQ, QRQ, QSQ, PQP, PRP* |
| 16 | (13) | (12) | (12) | (1) | (123) | (1) | (12) | (123) | (13) | (12) | *S, QR, QPQ, QSQ, PQP, PRP, PSP* |
| 17 | (13) | (12) | (1) | (12) | (123) | (12) | (12) | (13) | (123) | (1) | *R, QS, QPQ, QRQ, PQP, PRP, PSP* |
| 18 | (13) | (1) | (12) | (12) | (13) | (12) | (1) | (123) | (123) | (12) | *Q, RS, RPR, RQR, PQP, PRP, PSP* |
| 19 | (12) | (12) | (23) | (13) | (1) | (123) | (123) | (123) | (132) | (132) | *R, SPS, SQS, SRP, SRQ, SRSRS* |
| 20 | (12) | (23) | (12) | (13) | (123) | (132) | (132) | (1) | (132) | (123) | *Q, SPS, SRS, SQP, SQR, SQSQS* |
| 21 | (12) | (23) | (13) | (12) | (123) | (123) | (123) | (132) | (1) | (132) | *Q, RPR, RSR, RQP, RQS, RQRQR* |
| 22 | (23) | (12) | (12) | (13) | (132) | (1) | (132) | (132) | (123) | (132) | *P, RQ, QPS, QSQ, SQS, SRS* |
| 23 | (23) | (12) | (13) | (12) | (132) | (132) | (123) | (123) | (132) | (1) | *P, SQ, QPR, QRQ, RQR, RSR* |
| 24 | (23) | (13) | (12) | (12) | (123) | (123) | (1) | (132) | (132) | (123) | *P, SR, QPR, QRQ, QSQ, RQR* |
| 25 | (1) | (12) | (13) | (23) | (12) | (132) | (132) | (13) | (23) | (123) | *P, S, RPR, RQR, QSR, RSPSR, RSRSR* |
| 26 | (12) | (1) | (13) | (23) | (12) | (13) | (132) | (132) | (123) | (23) | *Q, S, RPR, RQR, PSR, RSQSR, RSRSR* |
| 27 | (12) | (13) | (1) | (23) | (132) | (13) | (23) | (12) | (123) | (132) | *R, S, QPQ, QRQ, PSQ, QSQSQ, QSRSQ* |
| 28 | (12) | (13) | (23) | (1) | (132) | (132) | (23) | (123) | (12) | (13) | *R, S, QPQ, QSQ, PRQ, QRQRQ, QRSRQ* |
| 29 | (13) | (12) | (12) | (12) | (123) | (1) | (1) | (123) | (123) | (1) | *QR, QS, QPQ, PQP, PRP, PSP* |
| 30 | (12) | (13) | (12) | (12) | (132) | (123) | (1) | (1) | (1) | (123) | *RP, SP, QPQ, QRQ, QSQ, PQP* |
| 31 | (12) | (12) | (13) | (12) | (1) | (132) | (123) | (132) | (1) | (1) | *QP, SP, RPR, RQR, RSR, PRP* |
| 32 | (12) | (12) | (12) | (13) | (1) | (1) | (132) | (1) | (132) | (132) | *QP, RP, SPS, SQS, SRS, PSP* |
| 33 | (12) | (12) | (13) | (13) | (1) | (132) | (1) | (132) | (132) | (132) | *RS, QP, RPR, RQR, PRP, PSP* |
| 34 | (12) | (13) | (12) | (13) | (132) | (123) | (132) | (1) | (132) | (1) | *QS, RP, QPQ, QRQ, PQP, PSP* |
| 35 | (12) | (13) | (13) | (12) | (132) | (1) | (123) | (132) | (1) | (123) | *QR, SP, QPQ, QSQ, PQP, PRP* |

Table 2. Colorings that will give rise to index 3 subgroups of *H*.

| no | *P* | *Q* | *R* | *S* | generators for the subgroup fixing color 1 | no | *P* | *Q* | *R* | *S* | generators for the subgroup fixing color 1 |
|---|---|---|---|---|---|---|---|---|---|---|---|
| 1 | (12)(34) | (13)(24) | (1) | (1) | *R, S, PRP, PSP, QPQP, QRQ, QSQ, PQRQP, PQSQP* | 19 | (12)(34) | (12)(34) | (14)(23) | (13)(24) | *QP, PRS, PSR, RPS, RQS* |
| 2 | (12)(34) | (1) | (13)(24) | (1) | *Q, S, PQP, PSP, RPRP, RQR, RSR, PRQRP, PRSRP* | 20 | (12)(34) | (14)(23) | (12)(34) | (13)(24) | *RP, PQS, PSQ, QPS, QRS* |
| 3 | (12)(34) | (1) | (1) | (13)(24) | *Q, R, PQP, PRP, SPSP, SQS, SRS, PSQSP, PSRSP* | 21 | (12)(34) | (14)(23) | (13)(24) | (12)(34) | *SP, PQR, PRQ, QPR, QSR* |
| 4 | (1) | (12)(34) | (13)(24) | (1) | *P, S, QPQ, QSQ, RPR, RQRQ, RSR, QRPRQ, QRSRQ* | 22 | (14)(23) | (12)(34) | (12)(34) | (13)(24) | *RQ, PQS, PRS, PSQ, QPS* |
| 5 | (1) | (12)(34) | (1) | (13)(24) | *P, R, QPQ, QRQ, SPS, SQSQ, SRS, QSPSQ, QSRSQ* | 23 | (14)(23) | (12)(34) | (13)(24) | (12)(34) | *SQ, PQR, PRQ, PSR, QPR* |
| 6 | (1) | (1) | (12)(34) | (13)(24) | *P, Q, RPR, RQR, SPS, SQS, SRSR, RSPSR, RSQSR* | 24 | (14)(23) | (13)(24) | (12)(34) | (12)(34) | *SR, PQR, PRQ, PSQ, QPR* |
| 7 | (12)(34) | (12)(34) | (1) | (13)(24) | *QP, R, PRP, SPSP, SQSP, SRS, PSRSP* | 25 | (1) | (12)(34) | (13)(24) | (14)(23) | *P, QPQ, QRS, QSR, RPR, RQS, SPS* |
| 8 | (12)(34) | (1) | (12)(34) | (13)(24) | *Q, RP, PQP, SPSP, SQS, SRSP, PSQSP* | 26 | (12)(34) | (1) | (13)(24) | (14)(23) | *Q, PQP, PRS, PSR, RPS, RQR, SQS* |
| 9 | (12)(34) | (1) | (13)(24) | (12)(34) | *Q, SP, PQP, RPRP, RQR, RSRP, PRQRP* | 27 | (12)(34) | (13)(24) | (1) | (14)(23) | *R, PQS, PRP, PSQ, QPS, QRQ, SRS* |
| 10 | (1) | (12)(34) | (12)(34) | (13)(24) | *P, RQ, QPQ, SPS, SQSQ, SRSQ, QSPSQ* | 28 | (12)(34) | (13)(24) | (14)(23) | (1) | *S, PQR, PRQ, PSP, QPR, QSQ, RSR* |
| 11 | (1) | (12)(34) | (13)(24) | (12)(34) | *P, SQ, QPQ, RPR, RQRQ, RSRQ, QRPRQ* | 29 | (13)(24) | (12)(34) | (12)(34) | (12)(34) | *RQ, SQ, PRQP, PSQP, QPQP* |
| 12 | (1) | (13)(24) | (12)(34) | (12)(34) | *P, SR, QPQ, QSRQ, RPR, RQRQ, QRPRQ* | 30 | (12)(34) | (13)(24) | (12)(34) | (12)(34) | *RP, SP, QPQP, QRQP, QSQP* |
| 13 | (12)(34) | (12)(34) | (13)(24) | (1) | *QP, S, PSP, RPRP, RQRP, RSR, PRSRP* | 31 | (12)(34) | (12)(34) | (13)(24) | (12)(34) | *QP, SP, RPRP, RQRP, RSRP* |
| 14 | (12)(34) | (13)(24) | (12)(34) | (1) | *RP, S, PSP, QPQP, QRQP, QSQ, PQSQP* | 32 | (12)(34) | (12)(34) | (12)(34) | (13)(24) | *QP, RP, SPSP, SQSP, SRSP* |
| 15 | (12)(34) | (13)(24) | (1) | (12)(34) | *R, SP, PRP, QPQP, QRQ, QSQP, PQRQP* | 33 | (12)(34) | (12)(34) | (13)(24) | (13)(24) | *QP, SR, PSRP, RPRP, RQRP* |
| 16 | (13)(24) | (12)(34) | (12)(34) | (1) | *RQ, S, PRQP, PSP, QPQP, QSQ, PQSQP* | 34 | (12)(34) | (13)(24) | (12)(34) | (13)(24) | *RP, SQ, PSQP, QPQP, QRQP* |
| 17 | (13)(24) | (12)(34) | (1) | (12)(34) | *R, SQ, PRP, PSQP, QPQP, QRQ, PQRQP* | 35 | (12)(34) | (13)(24) | (13)(24) | (12)(34) | *RQ, SP, PRQP, QPQP, QSQP* |
| 18 | (13)(24) | (1) | (12)(34) | (12)(34) | *Q, SR, PQP, PSRP, RPRP, RQR, PRQRP* | | | | | | |

Table 3. Colorings that will give rise to index 4 subgroups of *H*, where $\pi(H) \cong V$.

| no | *P* | *Q* | *R* | *S* | generators for the subgroup fixing color 1 | no | *P* | *Q* | *R* | *S* | generators for the subgroup fixing color 1 |
|---|---|---|---|---|---|---|---|---|---|---|---|
| 1 | (12) | (13)(24) | (1) | (1) | *R, S, QPQ, QRQ, QSQ, PRP, PSP, PQPQP, PQRQP, PQSQP* | 101 | (12) | (13)(24) | (12)(34) | (14)(23) | *SPS, SQP, SQR, SRQ, SQSRS, SRPRS* |
| 2 | (12) | (1) | (13)(24) | (1) | *Q, S, RPR, RQR, RSR, PQP, PSP, PRPRP, PRQRP, PRSRP* | 102 | (12) | (13)(24) | (14)(23) | (12)(34) | *RPR, RQP, RQS, QSR, RSPSR, RQRSR* |
| 3 | (12) | (1) | (1) | (13)(24) | *Q, PRP, PQP, SPS, SQS, SRS, PSPSP, PSQSP, PSRSP* | 103 | (12) | (13)(24) | (13)(24) | (13)(24) | *QR, QS, QPQ, PRQP, PSQP, PQPQP* |
| 4 | (1) | (12) | (13)(24) | (1) | *P, S, RPR, RQR, RSR, QPQ, QSQ, QRPRQ, QRQRQ, QRSRQ* | 104 | (12) | (13)(24) | (13)(24) | (14)(23) | *QR, QPQ, PSQ, SPS, QSQS, QSRS* |
| 5 | (1) | (12) | (1) | (13)(24) | *P, R, SPS, SQS, SRS, QPQ, QRQ, QSPSQ, QSQSQ, QSRSQ* | 105 | (12) | (13)(24) | (14)(23) | (13)(24) | *QS, QPQ, PRQ, RPR, QRQR, QRSR* |
| 6 | (1) | (1) | (12) | (13)(24) | *P, Q, RPR, RQR, SPS, SQS, SRS, RSPSR, RSQSR, RSRSR* | 106 | (12) | (14)(23) | (13)(24) | (13)(24) | *SR, QPQ, PRQ, RPR, QSRQ, QRQR* |
| 7 | (12) | (12) | (1) | (13)(24) | *R, QP, SPS, SQS, SRS, PRP, PSPSP, PSQSP, PSRSP* | 107 | (12)(34) | (1) | (1) | (13) | *Q, R, PQP, PRP, PSP, SQS, SRS, SPQPS, SPRPS, SPSPS* |
| 8 | (12) | (1) | (12) | (13)(24) | *Q, RP, SPS, SQS, SRS, PQP, PSPSP, PSQSP, PSRSP* | 108 | (12)(34) | (1) | (13) | (1) | *Q, S, PQP, PRP, PSP, RQR, RSR, RPQPR, RPRPR, RPSPR* |
| 9 | (12) | (1) | (13)(24) | (12) | *Q, SP, RPR, RQR, RSR, PQP, PRPRP, PRQRP, PRSRP* | 109 | (12)(34) | (13) | (1) | (1) | *R, S, PQP, PRP, PSP, QRQ, QSQ, QPQPQ, QPRPQ, QPSPQ* |
| 10 | (1) | (12) | (12) | (13)(24) | *P, QR, SPS, SQS, SRS, QPQ, QSPSQ, QSQSQ, QSRSQ* | 110 | (12)(34) | (12) | (13)(24) | (1) | *RQR, RSR, RPQPR, RPSPR, PRPR, QRPR, SRPR* |
| 11 | (1) | (12) | (13)(24) | (12) | *P, SQ, RPR, RQR, RSR, QPQ, QRPRQ, QRQRQ, QRSRQ* | 111 | (12)(34) | (12) | (1) | (13)(24) | *R, QP, SQS, QRQ, PRP, PSPS, QPQP, QPRPQ, QPSPQ* |
| 12 | (1) | (13)(24) | (12) | (12) | *P, SR, QPQ, QRQ, QSQ, RPR, RQPQR, RQRQR, RQSQR* | 112 | (12)(34) | (1) | (12) | (13)(24) | *Q, SQS, SRS, PSPS, RSPS, SPQPS, SPRPS, SPSQSPS* |
| 13 | (12) | (12) | (13)(24) | (1) | *S, QP, RPR, RQR, RSR, PSP, PRPRP, PRQRP, PRSRP* | 113 | (12)(34) | (13) | (24) | (1) | *R, S, PQP, PSP, QRQ, QSQ, PRPQ, PRQRP, PRSRP* |
| 14 | (12) | (13)(24) | (12) | (1) | *S, RP, QPQ, QRQ, QSQ, PSP, PQPQP, PQRQP, PQSQP* | 114 | (12)(34) | (13) | (1) | (24) | *R, S, PQP, PRP, QRQ, QSQ, PSPQ, PSQSP, PSRSP* |
| 15 | (12) | (13)(24) | (1) | (12) | *R, SP, QPQ, QRQ, QSQ, PRP, PQPQP, PQRQP, PQSQP* | 115 | (12)(34) | (1) | (13) | (24) | *Q, S, PQP, PRP, RQR, RSR, PSPR, PSQSP, PSRSP* |
| 16 | (13)(24) | (12) | (12) | (1) | *S, RQ, PQP, PRP, PSP, QSQ, QPQPQ, QPRPQ, QPSPQ* | 116 | (12)(34) | (13) | (12)(34) | (1) | *S, PR, PQP, PSP, QSQ, QRPQ, QPQPQ, QPSPQ* |
| 17 | (13)(24) | (12) | (1) | (12) | *R, SQ, PQP, PRP, QRQ, PSPQ, PSQSP, PSRSP* | 117 | (12)(34) | (13) | (1) | (12)(34) | *R, PS, PQP, PRP, QRQ, QSPQ, QPQPQ, QPRPQ* |
| 18 | (13)(24) | (1) | (12) | (12) | *Q, SR, PQP, PRP, RQR, PSP, PSQSP, PSRSP* | 118 | (12)(34) | (1) | (13) | (12)(34) | *Q, PS, PQP, PRP, RQR, RSPR, RPQPR, RPRPR* |
| 19 | (12) | (12) | (12)(34) | (13)(24) | *SPS, SQS, PSRS, QSRS, RSRS, SRPRS, SRQRS* | 119 | (12)(34) | (13) | (13)(24) | (1) | *S, RQ, PQP, PSP, QSQ, PRPQ, PRQRP, PRSRP* |
| 20 | (12) | (12)(34) | (12) | (13)(24) | *SPS, SRS, PSQS, QSQS, RSQS, SQPQS, SQRQS* | 120 | (12)(34) | (13) | (1) | (13)(24) | *PRP, SPSP, PSQSP, PSRSP, PSPQPSP, PSPRPSP, PQSPSP, PSPSRSPSP* |
| 21 | (12) | (12)(34) | (13)(24) | (12) | *RPR, RSR, PRQR, QRQR, SRQR, RQPQR, RQSQR* | 121 | (12)(34) | (1) | (13) | (13)(24) | *Q, SR, PQP, PRP, PSP, RQR, RPQPR, RPRPR, RPSPR* |
| 22 | (12)(34) | (12) | (12) | (13)(24) | *SQS, SRS, PSPS, QSPS, RSPS, SPQPS, SPRPS* | 122 | (12)(34) | (13) | (14)(23) | (1) | *S, PQP, PSP, QRP, RQR, RSR, PRPR, PRSRP* |
| 23 | (12)(34) | (12) | (13)(24) | (12) | *RS, RQR, RPRP, RPQPR, RPSPR, RPRQRPR, RPRSRPR* | 123 | (12)(34) | (13) | (1) | (14)(23) | *R, PQP, PRP, QSP, SQS, SRS, PSPS, PSRSP* |
| 24 | (12)(34) | (13)(24) | (12) | (12) | *QRQ, QSQ, PQPQ, RQPQ, SQPQ, QPRPQ, QPSPQ* | 124 | (12)(34) | (1) | (13) | (14)(23) | *Q, PQP, PRP, RSP, SQS, SRS, PSPS, PSQSP* |
| 25 | (1) | (12) | (34) | (13)(24) | *P, R, SPS, SQS, QSRS, SRPRS, SRQRS, SRSPSRS, SRSRSRS* | 125 | (12)(34) | (12)(34) | (13) | (1) | *S, PQ, PRP, PSP, RSR, RQPR, RPRPR, RPSPR* |
| 26 | (12) | (1) | (34) | (13)(24) | *Q, R, SPS, SQS, PSRS, SRPRS, SRQRS, SRSQSRS, SRSRSRS* | 126 | (12)(34) | (12)(34) | (1) | (13) | *R, PQ, PRP, PSP, SRS, SQPS, SPRPS, SPSPS* |
| 27 | (12) | (34) | (1) | (13)(24) | *Q, R, SPS, SRS, PSQS, SQPQS, SQRQS, SQSQSQS, SQSRSQS* | 127 | (12)(34) | (1) | (12)(34) | (13) | *Q, PR, PQP, PSP, SQS, SRPS, SPQPS, SPSPS* |
| 28 | (12) | (34) | (13)(24) | (1) | *Q, S, RPR, RSR, PRQR, RQPQR, RQSQR, RQRQRQR, RQRSRQR* | 128 | (12)(34) | (13)(24) | (12) | (1) | *S, QRQ, QSQ, QPQP, QPQR, QPRPQ, QPSPQ, QPQSQPQ* |
| 29 | (1) | (12) | (13)(24) | (34) | *P, S, RPR, RQR, QPQ, QSQ, RSRQ, RSPSR, RSQSR* | 129 | (12)(34) | (13)(24) | (1) | (12) | *R, QRQ, QSQ, PQPQ, SQPQ, QPRPQ, QPSPQ, QPQRQPQ* |
| 30 | (1) | (13)(24) | (12) | (34) | *P, S, QPQ, QRQ, RPR, RSR, QSQR, QSPSQ, QSRSQ* | 130 | (12)(34) | (1) | (13)(24) | (12) | *Q, RQR, RSR, PRPR, SRPR, RPQPR, RPSPR, RPRQRPR* |
| 31 | (1) | (13)(24) | (12) | (1) | *P, S, QPQ, QRQ, QSQ, RPR, RSR, RQPQR, RQRQR, RQSQR* | 131 | (12)(34) | (13)(24) | (13) | (1) | *S, RQ, PRP, PSP, QSQ, PQPQ, PQRQP, PQSQP* |
| 32 | (1) | (13)(24) | (1) | (12) | *P, R, QPQ, QRQ, QSQ, SPS, SRS, SQPQS, SQRQS, SQSQS* | 132 | (12)(34) | (13)(24) | (1) | (13) | *R, SQ, PRP, PSP, QRQ, PQPQ, PQRQP, PQSQP* |
| 33 | (1) | (1) | (13)(24) | (12) | *P, Q, RPR, RQR, RSR, SPS, SQS, SRPRS, SRQRS, SRSRS* | 133 | (12)(34) | (1) | (13)(24) | (13) | *PSP, RPRP, PRQRP, PRSRP, PQRPRP, PRPQPRP, PRPSPRP, PRPRSRPRP* |
| 34 | (1) | (12)(34) | (13)(24) | (12) | *P, RPR, RSR, QRQR, SRQR, RQPQR, RQSQR, RQRPRQR* | 134 | (12)(34) | (13)(24) | (14) | (1) | *S, QRQ, QSQ, QPR, PQPQ, QPSPQ, QPQRQPQ, QPQSQPQ* |
| 35 | (1) | (12)(34) | (12) | (13)(24) | *P, SPS, SRS, QSQS, RSQS, SQPQS, SQRQS, SQSPSQS* | 135 | (12)(34) | (13)(24) | (1) | (14) | *R, PRP, PSP, SQP, QRQ, QSQ, PQPQ, PQRQP* |
| 36 | (1) | (12) | (12)(34) | (13)(24) | *P, SPS, SQS, QSRS, RSRS, SRPRS, SRQRS, SRSPSRS* | 136 | (12)(34) | (1) | (13)(24) | (14) | *Q, PQP, PSP, SRP, RQR, RSR, PRPR, PRQRP* |
| 37 | (1) | (13)(24) | (12)(34) | (12) | *P, SR, QPQ, QSQ, RPR, QRQR, QRPRQ, QRSRQ* | 137 | (12)(34) | (12) | (34) | (13)(24) | *R, SQS, SRPS, PSPS, QSPS, SPQPS, SPSRSPS* |
| 38 | (1) | (13)(24) | (12) | (12)(34) | *P, SR, QPQ, QRQ, RPR, QSQR, QSPSQ, QSRSQ* | 138 | (12)(34) | (12) | (13)(24) | (34) | *S, RQR, RSPR, PRPR, QRPR, RPQPR, RPRSRPR* |
| 39 | (1) | (12) | (13)(24) | (12)(34) | *P, SQ, RPR, RQR, QPQ, RSRQ, RSPSR, RSQSR* | 139 | (12)(34) | (13)(24) | (12) | (34) | *S, QRQ, QPSQ, QPQP, QPQR, QPRPQ, QPQSQPQ* |
| 40 | (1) | (13)(24) | (13)(24) | (12) | *P, QR, QPQ, QSQ, SPS, SRQS, SQPQS, SQSQS* | 140 | (12)(34) | (12)(34) | (13)(24) | (12) | *RSR, RQPR, PRPR, QRPR, SRPR, RPSPR* |
| 41 | (1) | (13)(24) | (12) | (13)(24) | *P, QS, QPQ, QRQ, RPR, RSQR, RQPQR, RQRQR* | 141 | (12)(34) | (12)(34) | (12) | (13)(24) | *SRS, SQPS, PSPS, QSPS, RSPS, SPRPS* |
| 42 | (1) | (12) | (13)(24) | (13)(24) | *P, RS, RPR, RQR, QPQ, QSRQ, QRPRQ, QRQRQ* | 142 | (12)(34) | (12) | (12)(34) | (13)(24) | *SQS, SRPS, PSPS, QSPS, RSPS, SPQPS* |
| 43 | (1) | (13)(24) | (14)(23) | (12) | *P, QPQ, QSQ, SRQ, RPR, RSR, QRQR, QRPRQ* | 143 | (12)(34) | (13)(24) | (12)(34) | (12) | *QSQ, QRPQ, PQPQ, RQPQ, SQPQ, QPSPQ* |
| 44 | (1) | (13)(24) | (12) | (14)(23) | *P, QPQ, QRQ, RSQ, SPS, SRS, QSQS, QSPSQ* | 144 | (12)(34) | (13)(24) | (12) | (12)(34) | *QPSQ, QSRSQ, QSPRPSQ, QSQPQPSQ, QSQRQPSQ, QSQSQPSQ* |
| 45 | (1) | (12) | (13)(24) | (14)(23) | *P, SPS, SQS, SRQ, RSRS, SRPRS, SRSPSRS, SRSQSRS* | 145 | (12)(34) | (12) | (13)(24) | (12)(34) | *RQR, RSPR, RPQPR, PRPR, QRPR, SRPR* |
| 46 | (12) | (12)(34) | (13)(24) | (1) | *S, RPR, RSR, PRQR, QRQR, RQPQR, RQSQR, RQRSRQR* | 146 | (12)(34) | (13)(24) | (13)(24) | (12) | *QR, QSQ, PQPQ, SQPQ, QPSPQ, QPRQPQ* |
| 47 | (12) | (12)(34) | (1) | (13)(24) | *R, SPS, SRS, PSQS, QSQS, SQPQS, SQRQS, SQSRSQS* | 147 | (12)(34) | (13)(24) | (12) | (13)(24) | *QS, QRQ, PQPQ, RQPQ, QPRPQ, QPSQPQ* |
| 48 | (12) | (1) | (12)(34) | (13)(24) | *Q, SPS, SQS, PSRS, RSRS, SRPRS, SRQRS, SRSQSRS* | 148 | (12)(34) | (12) | (13)(24) | (13)(24) | *RS, RQR, PRPR, QRPR, RPQPR, RPSRPR* |
| 49 | (12) | (13)(24) | (34) | (1) | *R, S, QPQ, QSQ, PRP, PSP, QRQP, QRPRQ, QRSRQ* | 149 | (12)(34) | (13)(24) | (14)(23) | (12) | *RSR, RPQ, RQP, RQS, RPRQR, RPSPR* |
| 50 | (12) | (13)(24) | (1) | (34) | *R, S, QPQ, QRQ, PRP, PSP, QSQP, QSPSQ, QSRSQ* | 150 | (12)(34) | (13)(24) | (12) | (14)(23) | *SRS, SPQ, SQP, SQR, SPRPS, SPSQS* |
| 51 | (12) | (1) | (13)(24) | (34) | *Q, S, RPR, RQR, PQP, PSP, RSRP, RSPSR, RSQSR* | 151 | (12)(34) | (12) | (13)(24) | (14)(23) | *SQS, SPR, SRP, SRQ, SPQPS, SPSRS* |
| 52 | (12) | (13)(24) | (12)(34) | (1) | *S, RP, QPQ, QSQ, PSP, QRQP, QRPRQ, QRSRQ* | 152 | (12)(34) | (13) | (13) | (13) | *RQ, SQ, PQP, PRP, PSPQ, PSQSP, PSRSP* |
| 53 | (12) | (13)(24) | (1) | (12)(34) | *R, SP, QPQ, QRQ, PRP, QSQP, QSPSQ, QSRSQ* | 153 | (12)(34) | (13) | (13) | (24) | *S, RQ, PQP, PRP, QSQ, PSPQ, PSQSP, PSRSP* |
| 54 | (12) | (1) | (13)(24) | (12)(34) | *Q, SP, RPR, RQR, PQP, RSRP, RSPSR, RSQSR* | 154 | (12)(34) | (13) | (24) | (13) | *R, SQ, PQP, PSP, QRQ, PRPQ, PRQRP, PRSRP* |
| 55 | (12) | (13)(24) | (13)(24) | (1) | *S, QR, QPQ, QSQ, PSP, PRQP, PQPQP, PQSQP* | 155 | (12)(34) | (24) | (13) | (13) | *Q, SR, PRP, PSP, RPQP, RQR, PQRQP, PQSQP* |
| 56 | (12) | (13)(24) | (1) | (13)(24) | *R, QS, QPQ, QRQ, PRP, PSQP, PQPQP, PQRQP* | 156 | (12)(34) | (13) | (13) | (12)(34) | *PS, RQ, PQP, PRP, QSPQ, QPQPQ, QPRPQ* |
| 57 | (12) | (1) | (13)(24) | (13)(24) | *Q, RS, RPR, RQR, PQP, PSRP, PRPRP, PRQRP* | 157 | (12)(34) | (13) | (12)(34) | (13) | *PR, SQ, PQP, PSP, QRPQ, QPQPQ, QPSPQ* |
| 58 | (12) | (13)(24) | (14)(23) | (1) | *S, QPQ, QSQ, PRQ, RPR, RSR, QRQR, QRSRQ* | 158 | (12)(34) | (12)(34) | (13) | (13) | *SR, PQ, PRP, PSP, RQPR, RPRPR, RPSPR* |
| 59 | (12) | (13)(24) | (1) | (14)(23) | *R, QPQ, QRQ, PSQ, SPS, SRS, QSQS, QSRSQ* | 159 | (12)(34) | (13) | (13) | (13)(24) | *RQ, SQ, PQP, PRP, PSP, QPQPQ, QPRPQ, QPSPQ* |
| 60 | (12) | (1) | (13)(24) | (14)(23) | *Q, RPR, RQR, PSR, SPS, SQS, RSRS, RSQSR* | 160 | (12)(34) | (13) | (13)(24) | (13) | *RQ, SQ, PQP, PSP, PRPQ, PRQRP, PRSRP* |
| 61 | (12) | (12) | (12) | (13)(24) | *QP, RP, SPS, SQS, SRS, PSPSP, PSQSP, PSRSP* | 161 | (12)(34) | (13)(24) | (13) | (13) | *RQ, SQ, PRP, PSP, PQPQ, PQRQP, PQSQP* |
| 62 | (12) | (12) | (13)(24) | (12) | *QP, SP, RPR, RQR, RSR, PRPRP, PRQRP, PRSRP* | 162 | (12)(34) | (13) | (13) | (14)(23) | *PQP, PRP, QSP, RSP, SQS, SRS, PSPS* |
| 63 | (12) | (13)(24) | (12) | (12) | *RP, SP, QPQ, QRQ, QSQ, PQPQP, PQRQP, PQSQP* | 163 | (12)(34) | (13) | (14)(23) | (13) | *PQP, PSP, QRP, SRP, RQR, RSR, PRPR* |
| 64 | (12) | (12) | (34) | (13)(24) | *R, SPS, SQS, PSRS, QSRS, SRPRS, SRQRS, SRSRSRS* | 164 | (12)(34) | (14)(23) | (13) | (13) | *PRP, PSP, RQP, SQP, QRQ, QSQ, PQPQ* |
| 65 | (12) | (12) | (13)(24) | (34) | *S, QP, RPR, RQR, PSP, RSRP, RSPSR, RSQSR* | 165 | (12)(34) | (13) | (24) | (12)(34) | *R, PS, PQP, QRQ, PRPQ, PRSQ, PRQRP* |
| 66 | (12) | (13)(24) | (12) | (34) | *S, RP, QPQ, QRQ, PSP, QSQP, QSPSQ, QSRSQ* | 166 | (12)(34) | (13) | (12)(34) | (24) | *S, PR, PQP, QSQ, PSPQ, PSRQ, PSQSP* |
| 67 | (12) | (13)(24) | (12)(34) | (12) | *RP, SP, QPQ, QSQ, QRQP, QRPRQ, QRSRQ* | 167 | (12)(34) | (12)(34) | (13) | (24) | *S, PQ, PRP, RSR, PSPR, PSQR, PSRSP* |
| 68 | (12) | (13)(24) | (12) | (12)(34) | *RP, SP, QPQ, QRQ, QSQP, QSPSQ, QSRSQ* | 168 | (12)(34) | (13) | (24) | (13)(24) | *R, SQ, PQP, QRP, QRQ, PSRP, PRPQ* |
| 69 | (12) | (12) | (13)(24) | (12)(34) | *QP, SP, RPR, RQR, RSRP, RSPSR, RSQSR* | 169 | (12)(34) | (13) | (13)(24) | (24) | *S, RQ, PQP, QSQ, PSRP, PRPQ, PRQRP* |
| 70 | (12) | (13)(24) | (13)(24) | (12) | *QR, SP, QPQ, QSQ, PRQP, PQPQP, PQSQP* | 170 | (12)(34) | (13)(24) | (13) | (24) | *S, RQ, PRP, QSQ, PSQP, PQPQ, PQRQP* |
| 71 | (12) | (13)(24) | (12) | (13)(24) | *QS, RP, QPQ, QRQ, PSQP, PQPQP, PQRQP* | 171 | (12)(34) | (13) | (24) | (14)(23) | *R, PQP, QSP, SRP, PRPSP, PRQRP, PSRSP* |
| 72 | (12) | (12) | (13)(24) | (13)(24) | *RS, QP, RPR, RQR, PSRP, PRPRP, PRQRP* | 172 | (12)(34) | (13) | (14)(23) | (24) | *S, PQP, QRP, RSP, PRPSP, PRSRP, PSQSP* |
| 73 | (12) | (13)(24) | (14)(23) | (12) | *QPQ, QSQ, PRQ, SRQ, RPR, RSR, QRQR* | 173 | (12)(34) | (14)(23) | (13) | (24) | *S, QRQ, QPR, PSQ, QPQSQ, QPSPQ, QSRSQ* |
| 74 | (12) | (13)(24) | (12) | (14)(23) | *QPQ, QRQ, PSQ, RSQ, SPS, SRS, QSQS* | 174 | (12)(34) | (12)(34) | (12)(34) | (13) | *PQ, PR, PSP, SQPS, SRPS, SPSPS* |
| 75 | (12) | (12) | (13)(24) | (14)(23) | *RPR, RQR, PSR, QSR, SPS, SQS, RSRS* | 175 | (12)(34) | (12)(34) | (13) | (12)(34) | *PS, PQ, PRP, RQPR, RSPR, RPRPR* |
| 76 | (12) | (34) | (12) | (13)(24) | *Q, SPS, SRS, PSQS, RSQS, SQPQS, SQRQS, SQSQSQS* | 176 | (12)(34) | (13) | (12)(34) | (12)(34) | *PR, PS, PQP, QRPQ, QSPQ, QPQPQ* |
| 77 | (12) | (34) | (13)(24) | (12) | *Q, RPR, RSR, PRQR, SRQR, RQPQR, RQSQR, RQRQRQR* | 177 | (12)(34) | (12)(34) | (13)(24) | (13) | *SR, PQ, PSP, PRPR, PRQR, PRSRP* |
| 78 | (12) | (13)(24) | (34) | (12) | *R, SP, QPQ, QSQ, PRP, QRQP, QRPRQ, QRSRQ* | 178 | (12)(34) | (12)(34) | (13) | (13)(24) | *SR, PQ, PRP, PSPR, PSQR, PSRSP* |
| 79 | (12) | (34) | (34) | (13)(24) | *Q, R, SPS, SRQS, PSQS, SQPQS, SQSQSQS, SQSRSQS* | 179 | (12)(34) | (13) | (12)(34) | (13)(24) | *PR, SQ, PQP, PSPQ, PSRQ, PSQSP* |
| 80 | (12) | (34) | (13)(24) | (34) | *Q, S, RPR, RSQR, PRQR, RQPQR, RQRQRQR, RQRSRQR* | 180 | (12)(34) | (12)(34) | (14)(23) | (13) | *PQ, PSP, SRP, RSR, PRPR, PRQR* |
| 81 | (12) | (13)(24) | (34) | (34) | *R, S, QPQ, PRP, PSP, QSRQ, QRQP, QRPRQ* | 181 | (12)(34) | (12)(34) | (13) | (14)(23) | *PQ, PRP, RSP, SRS, PSPS, PSQS* |
| 82 | (12) | (12)(34) | (13)(24) | (34) | *S, RPR, RSQR, PRQR, QRQR, RQPQR, RQRSRQR* | 182 | (12)(34) | (13) | (12)(34) | (14)(23) | *PR, PQP, QSP, SQS, PSPS, PSRS* |
| 83 | (12) | (12)(34) | (34) | (13)(24) | *R, SPS, SRQS, PSQS, QSQS, SQPQS, SQSRSQS* | 183 | (12)(34) | (13)(24) | (12)(34) | (13) | *SQ, PR, PSP, PQPQ, PQRQ, PQSQP* |
| 84 | (12) | (34) | (12)(34) | (13)(24) | *Q, SPS, SRQS, PSQS, RSQS, SQPQS, SQSQSQS* | 184 | (12)(34) | (13)(24) | (13) | (12)(34) | *PS, RQ, PRP, PQPQ, PQSQ, PQRQP* |
| 85 | (12) | (13)(24) | (12)(34) | (34) | *S, RP, QPQ, PSP, QSRQ, QRQP, QRPRQ* | 185 | (12)(34) | (13) | (13)(24) | (12)(34) | *PS, RQ, PQP, PRPQ, PRSQ, PRQRP* |
| 86 | (12) | (13)(24) | (34) | (12)(34) | *R, SP, QPQ, PRP, QSRQ, QRQP, QRPRQ* | 186 | (12)(34) | (13)(24) | (13)(24) | (13) | *RQ, SQ, PSP, PRQP, PQPQ, PQSQP* |
| 87 | (12) | (34) | (13)(24) | (12)(34) | *Q, RPR, RSQR, PRQR, SRQR, RQPQR, RQRQRQR* | 187 | (12)(34) | (13)(24) | (13) | (13)(24) | *RQ, SQ, PRP, PSQP, PQPQ, PQRQP* |
| 88 | (12) | (13)(24) | (13)(24) | (34) | *S, QR, QPQ, PSP, QSQP, QSRP, QSPSQ* | 188 | (12)(34) | (13) | (13)(24) | (13)(24) | *RQ, SQ, PQP, PSRP, PRPQ, PRQRP* |
| 89 | (12) | (13)(24) | (34) | (13)(24) | *R, QS, QPQ, PRP, QRQP, QRSP, QRPRQ* | 189 | (12)(34) | (13)(24) | (14)(23) | (13) | *PSP, QRP, RQP, SRP, PQPRP, PQSQP* |
| 90 | (12) | (34) | (13)(24) | (13)(24) | *Q, RS, RPR, PRQR, RQPQR, RQSRQR, RQRQRQR* | 190 | (12)(34) | (13)(24) | (13) | (14)(23) | *PRP, QSP, RSP, SQP, PQPSP, PQRQP* |
| 91 | (12) | (13)(24) | (14)(23) | (34) | *S, RPR, RQP, QSR, RQRSR, RQSQR, RSPSR* | 191 | (12)(34) | (13) | (13)(24) | (14)(23) | *SQS, SPQ, SPR, SRP, SPSRS, SRQRS* |
| 92 | (12) | (13)(24) | (34) | (14)(23) | *R, QPQ, PSQ, SRQ, QRPRQ, QRQSQ, QSRSQ* | 192 | (12)(34) | (14)(23) | (12)(34) | (13) | *PR, PSP, SQP, QSQ, PQPQ, PQRQ* |
| 93 | (12) | (34) | (13)(24) | (14)(23) | *Q, RPR, RQS, PSR, RQPQR, RQRSR, RSQSR* | 193 | (12)(34) | (14)(23) | (13) | (12)(34) | *PS, PRP, RQP, QRQ, PQPQ, PQSQ* |
| 94 | (12) | (12)(34) | (12)(34) | (13)(24) | *SPS, SRQS, PSQS, QSQS, RSQS, SQPQS* | 194 | (12)(34) | (13) | (14)(23) | (12)(34) | *PS, PQP, QRP, RQR, PRSR, PRPRPP* |
| 95 | (12) | (12)(34) | (13)(24) | (12)(34) | *RPR, RSQR, PRQR, QRQR, SRQR, RQPQR* | 195 | (12)(34) | (14)(23) | (13)(24) | (13) | *PSP, QRP, RQP, SQP, PQPRP, PRSRP* |
| 96 | (12) | (13)(24) | (12)(34) | (12)(34) | *RP, SP, QPQ, QSRQ, QRQP, QRPRQ* | 196 | (12)(34) | (14)(23) | (13) | (13)(24) | *PRP, QSP, RQP, SQP, PQPSP, PSRSP* |
| 97 | (12) | (13)(24) | (13)(24) | (12)(34) | QR, SP, QPQ, QSQP, QSRP, QSPSQ | 197 | (12)(34) | (13) | (14)(23) | (13)(24) | RQR, RPQ, RPRSR, RPS, PSR, RSQSR |
| 98 | (12) | (13)(24) | (12)(34) | (13)(24) | QS, RP, QPQ, QRQP, QRSP, QRPRQ | 198 | (12)(34) | (14)(23) | (14)(23) | (13) | RQ, PSP, SQP, QSQ, PRQP, PQPQ |
| 99 | (12) | (12)(34) | (13)(24) | (13)(24) | RS, SPS, SQPQS, SQSRQS, SRPRQS, SRQRQS | 199 | (12)(34) | (14)(23) | (13) | (14)(23) | SQ, PRP, RQP, QRQ, PSQP, PQPQ |
| 100 | (12) | (12)(34) | (13)(24) | (14)(23) | SPS, SQR, SRP, SRQ, SQPQS, SQSRS | 200 | (12)(34) | (13) | (14)(23) | (14)(23) | SR, PQP, QRP, RQR, PSRP, PRPR |

Table 4. Colorings that will give rise to index 4 subgroups of $H$, where $\pi(H) \cong D_4$.

| no | $P$ | $Q$ | $R$ | $S$ | generators for the subgroup fixing color 1 | no | $P$ | $Q$ | $R$ | $S$ | generators for the subgroup fixing color 1 |
|---|---|---|---|---|---|---|---|---|---|---|---|
| 1 | (1) | (12) | (13) | (12)(34) | P, PSP, PQPQP, PQRQP, PQSQP, PRPRP, PRQRP, PRSPSRP, PRSQSRP, PRSRSRP | 122 | (12) | (34) | (12) | (13) | QR, RPR, RSPSR, RSRSR, RQPQR, RQSQR, RSQPQSR, RSQRQSR, RSQSQSR |
| 2 | (1) | (12) | (12)(34) | (13) | P, RPR, RQR, RSR, QPS, QRQ, QSQ, SQS, SRS | 123 | (12) | (12) | (34) | (13) | SRP, PSP, PRPRP, PRQRP, PRSQSRP, PRSRSRP, PQPQP, PQRQP, PQSQP |
| 3 | (1) | (12)(34) | (12) | (13) | P, RPR, RQR, RSR, QPQ, QRS, QSQ, SPS, SQS | 124 | (12) | (13) | (14) | (12) | SR, SQS, SPRPS, SRPRS, SRQRS, SPSQPS, SPQPQPS, SPQRQPS |
| 4 | (1) | (12) | (13) | (13)(24) | P, PQP, PRP, PSPSP, PSRSP, PSQPQSP, PSQSRQSP, PSQRPRQSP, PSQRQRQSP | 125 | (12) | (13) | (12) | (14) | RP, SPS, SQS, SRS, PQP, PSP, QPQ, QRQ, QSQ |
| 5 | (1) | (12) | (13)(24) | (13) | P, PQP, PSP, PRPRP, PRQRP, PRSQSRP, PRSRPSRP, PRSPQPSRP, PRSPSPSRP | 126 | (12) | (12) | (13) | (14) | QP, SPS, SQS, SRS, PRP, PSP, RPR, RQR, RSR |
| 6 | (1) | (13)(24) | (12) | (13) | P, PQP, PSP, PRQRP, PRSPRP, PRPRQPRP, PRPQPQPRP, PRPQSQPRP | 127 | (12) | (13) | (12) | (23) | RP, SPS, SQS, SRS, PQP, PSQ, QPQ, QRQ |
| 7 | (1) | (12) | (13) | (14)(23) | P, SPS, SQS, SRS, QPR, QRQ, QSR, RQR | 128 | (12) | (12)(34) | (13)(24) | (13) | SQS, QPS, SRQPS, SPSPS, SPQPQPS, SPRPRPS, SPRQRPS, SPRSRPS |
| 8 | (1) | (12) | (14)(23) | (13) | P, PQP, PSP, PRPRP, PRQRP, PRSPSRP, PRSRSRP, PRSQPQSRP, PRSQRQSRP, PRSQSQSRP | 129 | (12) | (12)(34) | (13) | (13)(24) | PSP, PRPRP, PRSRP, PRQS, PRQPQRP, PRQRQRP, PRQSQP, PRQSRSQRP |
| 9 | (1) | (14)(23) | (12) | (13) | P, RPR, RQR, RSR, QPS, QRS, QSQ, SQS | 130 | (12) | (12)(34) | (14)(23) | (13) | SPS, RQS, SQPQS, SQSQS, SQRPRS, SQRQRS, SRSRS |
| 10 | (1) | (12) | (13) | (24) | P, PQP, PRP, PSPSP, PSRSP, PSQPQSP, PSQSQSP, PSQRPRQSP, PSQRQRQSP, PSQRSRQSP | 131 | (12) | (12)(34) | (13) | (14)(23) | RPR, SQR, RQPQR, RQRQR, RQSPSR, RQSQSQR, RSQSR, RSRSR |
| 11 | (1) | (13) | (24) | (12) | P, PQP, PRPRP, PRSRP, PSPQRP, PSQSP, PSRQRP, PRQSQRP | 132 | (12) | (13) | (12)(34) | (14)(23) | QSR, QRQ, QPQPQ, QPRPQ, QPSPQ, QSPSQ, QSQSQ, QSRPRSQ, QSRSRSQ |
| 12 | (1) | (24) | (12) | (13) | P, PQP, PRP, PSPSP, PSQSP, PSRPRSP, PSRSQRSP, PSRQPQRSP, PSRQRQRSP | 133 | (12) | (13)(24) | (12)(34) | (13) | QPQ, PQRQ, QRPQRQ, QRPRQ, QRSRQ, QRQRPQRQ, QRQSQRQ |
| 13 | (1) | (12) | (34) | (23) | P, R, S, QRPQ, QPQSQ, QPSPQ, QSPSQ, QSRSQ | 134 | (12) | (13)(24) | (13) | (12)(34) | RS, SPS, SRQRS, SQPQS, SQRQS, SQSPRS, SQSQSQS, SQSRSQS |
| 14 | (1) | (12) | (23) | (34) | P, R, S, QPQ, QSQ, QRQPRQ, QRSPRQ, QRPRPRQ | 135 | (12) | (13) | (13)(24) | (12)(34) | QRQ, QSQ, QPRS, QPQPQ, QPSPQ, QPRPRPQ, QPRQRPQ, QPRSPSRPQ, QPRSRSRPQ |
| 15 | (1) | (23) | (12) | (13) | P, SPS, SQS, SRS, QPQ, QRQ, QSQ, RPR, RQR, RSR | 136 | (12) | (14)(23) | (12)(34) | (13) | RP, QPQ, QSQ, SRQ, PSP, QRQP, QRPRQ |
| 16 | (1) | (12) | (23) | (24) | P, R, S, QPQ, QRPSQ, QRQSQ, QRSRQ, QSRSQ | 137 | (12) | (14)(23) | (13) | (12)(34) | PS, SRS, SPS, SQPS, SPRPS, SQPQS, SQRQS |
| 17 | (12) | (13) | (12)(34) | (1) | S, SPS, SQS, SRPRS, SRQRS, SRSQSRS, SRSRSRS, SRSPQPSRS, SRSPRPSRS, SRSPSPSRS | 138 | (12) | (13) | (14)(23) | (12)(34) | SP, RPR, RQR, QSR, PQP, RSRP, RSPSR |
| 18 | (12) | (13) | (1) | (12)(34) | R, SP, PQP, PRP, QPQ, QRQ, QSRSQ, QSPSQ, QSQSQ | 139 | (12) | (13)(24) | (13)(24) | (13) | RQ, RPR, RQSQR, RQPQR, RSQSR, RSRPSR, RSPQPSR, RSPSPSR |
| 19 | (12) | (1) | (13) | (12)(34) | Q, SP, PQP, PRP, RPR, RQR, RSQSR, RSPSR, RSRSR | 140 | (12) | (13)(24) | (13) | (13)(24) | RQ, SQ, QPQ, PRP, PQSP, PQPQP, PQRQP |
| 20 | (12) | (13) | (13)(24) | (1) | S, SPS, SQS, SRSRS, SRPQPRS, SRPRPRS, SRPSPRS, SRQPQRS, SRQRQRS, SRQSQRS | 141 | (12) | (13) | (13)(24) | (13)(24) | RQ, SQ, QPQ, PQP, PRSP, PRPRP, PRQRP |
| 21 | (12) | (13) | (1) | (13)(24) | R, PQP, PRP, SQS, QPQ, QRQ, PSRSP, PSPSP, PSQSP | 142 | (12) | (14)(23) | (14)(23) | (13) | RQ, RPR, RQPQR, RQSQR, RSQPSR, RSRSR, RSPRPSR, RSPSPSR |
| 22 | (12) | (1) | (13) | (13)(24) | Q, QPQ, QRQ, QSQSQ, QSRSQ, QSPQPSQ, QSPSRPSQ, QSPRPRPSQ, QSPRQRPSQ | 143 | (12) | (14)(23) | (13) | (14)(23) | QS, QPQ, QRQ, PQR, PRP, PSP, RPR, RSR |
| 23 | (12) | (13) | (14)(23) | (1) | S, SPS, SQS, SRQRS, SRSRS, SRPRQPRS, SRPSPRS, SRPQPQPRS, SRPQSQPRS | 144 | (12) | (13) | (14)(23) | (14)(23) | RS, RPR, RQR, PQP, PRQ, PSQ, QPQ |
| 24 | (12) | (13) | (1) | (14)(23) | R, SPS, SQS, SRS, PQP, PRP, PSQ, QPQ, QRQ | 145 | (12) | (13)(24) | (14)(23) | (13) | RPR, RSR, RQP, RQSQR, QRQR, SRQR, RQRPRQR |
| 25 | (12) | (1) | (13) | (14)(23) | Q, SPS, SQS, SRS, PQR, PRP, PSR, RPR | 146 | (12) | (13)(24) | (13) | (14)(23) | RQ, SPS, SRS, SQP, QPQ, QSQ, SQRQS, SQSQS |
| 26 | (12) | (12)(34) | (13) | (1) | S, SPS, SQS, SRQRS, SRSRS, SRPRQPRS, SRPSQPRS, SRPQPQPRS | 147 | (12) | (13) | (13)(24) | (14)(23) | SPS, SQS, SRP, SRQRS, QSRS, RSRS, SRSPSRS |
| 27 | (12) | (12)(34) | (1) | (13) | RPR, RQR, RSR, SP, PQP, PRQ, QPQ, QSQ | 148 | (12) | (14)(23) | (13) | (13)(24) | SR, QPQ, QRQ, PSQ, RPR, QSQR, QSRSQ |
| 28 | (12) | (1) | (12)(34) | (13) | Q, QPQ, QSQ, QRQPRQ, QRSRQ, QRPSRPRQ, QRPRPRPRQ, QRPRQRPRQ | 149 | (12) | (13) | (14)(23) | (13)(24) | SQ, RPR, RQR, PSR, QPQ, RSRQ, RSQSR |
| 29 | (12) | (13)(24) | (13) | (1) | S, SRS, SQPQS, SQRQS, SQSQS, SPRQPS, SPSPS, SPQPQPS, SPQSQPS | 150 | (23) | (12)(34) | (13)(24) | (24) | P, S, QPR, QSRQ, RQRQ, RSR, QRPRQ |
| 30 | (12) | (13)(24) | (1) | (13) | R, RSR, RQPQR, RQSQR, RPQPR, RPRPR, RPSPR, RQRPRQR, RQRQRQR, RQRSRQR | 151 | (12) | (13) | (13) | (24) | S, RQ, PQP, PRP, QPQ, QSQ, PSPSP, PSQSP, PSRSP |
| 31 | (12) | (1) | (13)(24) | (13) | Q, QSQ, QRPRQ, QRQRQ, QRSRQ, QPQPQ, QPSRPQ, QPRPRPQ, QPRQRPQ | 152 | (12) | (13) | (24) | (13) | SQ, SQPQS, SQRQS, SPQRS, SPRPS, SPSRS, SRPRS |
| 32 | (12) | (13)(24) | (14) | (1) | S, QPQ, QRQ, QSQ, PQR, PRP, PSP, RPR, RSR | 153 | (12) | (24) | (13) | (13) | Q, QPQ, QSRQ, QRQRQ, QRPSPRQ, QRPQPQPRQ, QRPQSQPRQ, |
| | | | | | | | | | | | QRPRQPRQ |
| 33 | (12) | (13)(24) | (1) | (14) | R, SP, QPQ, QSQ, PQP, PRP, QRPRQ, QRQRQ, QRSRQ | 154 | (12) | (13) | (13) | (34) | S, QR, QPQ, PQP, PRP, PSP, QSPSQ, QSQSQ, QSRSQ |
| 34 | (12) | (1) | (13)(24) | (14) | Q, RPR, RQR, RSR, PQS, PRS, PSP, SPS | 155 | (12) | (13) | (34) | (13) | RS, RPR, RQRP, RQRQ, RQPQR, RQSQR, RQRSRQR |
| 35 | (12) | (12) | (13) | (12)(34) | QP, SP, PRP, RPR, RQR, RSPSR, RSQSR, RSRSR | 156 | (12) | (34) | (13) | (13) | RS, SPS, SQPQS, SQRQS, SQSQS, SRQPRS, SRPRPRS, SRPSPRS |
| 36 | (12) | (12) | (12)(34) | (13) | RP, SPS, SQS, SRQ, PQP, PSP, SRPRS, SRSRS | 157 | (12) | (13) | (14) | (23) | S, SQS, SPSPSPS, SPSQSPS, SPSRSPS, SPRPS, SRQPS, SPQPQPS, SPQSQPS |
| 37 | (12) | (12)(34) | (12) | (13) | SPS, RS, SQRPRQS, SQRQRQS, SQRSRQS, SQPQS, SQSQS, SRPRS, SRQRS | 158 | (12) | (23) | (13) | (14) | Q, SPS, SQS, SRS, PQR, PRP, PSP, RPR, RSR |
| 38 | (12) | (12) | (13) | (13)(24) | PQ, SR, PRP, RPR, RQR, PSPSP, PSQSP, PSRSP | 159 | (12) | (13) | (24) | (14) | R, QPQ, QRQ, QSQ, SPS, SQS, SRP, SRQRS, SRSRS |
| 39 | (12) | (12) | (13)(24) | (13) | RS, RPR, RQR, PQP, PRP, PSQ, QPQ, QRQ | 160 | (12) | (13) | (14) | (34) | S, SPS, SQPQS, SQRQS, SRSQS, SQSQSQS, SQSPQPSQS, SQSPRPSQS, SQSPSPSQS |
| 40 | (12) | (13)(24) | (12) | (13) | PR, RSR, RQRPRQR, RQRQRQR, RQRSRQR, RQPQR, RQSQR, RPQPR, RPSPR | 161 | (12) | (13) | (34) | (14) | SPS, RQS, SRQRS, SRSRS, SQPQS, SQSQS, SQRPRS, SQRQRQS |
| 41 | (12) | (12) | (13) | (14)(23) | QP, SPS, SQS, SRS, PRP, PSR, RPR, RQR | 162 | (12) | (34) | (13) | (14) | Q, QSQ, QRPRQ, QPSRQ, QRSQSRQ, QRSRSRQ, QRQPQRQ, QRQRQRQ, QRQSQRQ |
| 42 | (12) | (12) | (14)(23) | (13) | SPS, SQS, SRQ, PSRS, RSRS, SRPRS, SRSQSRS | 163 | (12) | (13) | (23) | (34) | R, S, PQ, PSPSP, PSRSP, PSQPQSP, PSQRP, PSQSQSP |
| 43 | (12) | (14)(23) | (12) | (13) | RP, SP, QPQ, QSQ, PQP, QRPRQ, QRQRQ, QRSRQ | 164 | (12) | (24) | (34) | (13) | SRQ, QSQ, QPQPQ, QPRPQ, QPSPQ, QRPRQ, QRQRQ, QRSPSRQ, QRSRSRQ |
| 44 | (12) | (13) | (12) | (12)(34) | RP, SP, QPQ, QRQ, PQP, QSPSQ, QSQSQ, QSRSQ | 165 | (12) | (24) | (13) | (34) | SPS, SRS, SQPQS, SQR, PSQS, QSQS, SQSRSQS |
| 45 | (12) | (13) | (12)(34) | (12) | PQ, QRQ, QSRPRSQ, QSRQRSQ, QSRSRSQ, QSPSQ, QSQSQ, QPRPQ, QPSPQ | 166 | (12) | (34) | (34) | (13) | Q, R, SPS, PQP, PRP, PSP, SQPQS, SQRS, SQSQS |
| 46 | (12) | (12)(34) | (13) | (12) | SP, RPR, RQR, RSR, PQP, PRP, QPQ, QRQ, QSQ | 167 | (12) | (34) | (13) | (34) | QS, SPS, SRS, SQPQS, SQRQRQS, SQRSPRQS, SQRPQPRQS, SQRPRPRQS |
| 47 | (12) | (13) | (12) | (13)(24) | PR, RQR, RSPSR, RSQSR, RSRSR, RPSPR, RPQPQPR, RPQRQPR, RPQSQPR | 168 | (12) | (13) | (34) | (34) | R, SP, PQP, PRP, QPQ, QRSQ, QRPRQ, QRQRQ |
| 48 | (12) | (13) | (13)(24) | (12) | SQP, SPS, SRS, SQRPRQS, SQRQRQS, SQRSRQS, SQSQS, SQPQPQS, SQPRPQS | 169 | (12) | (34) | (34) | (23) | Q, R, S, PQP, PRP, PSQPSP, PSRPSP, PSPSPSP |
| 49 | (12) | (13)(24) | (13) | (12) | PS, RQ, PRP, QPQ, QSQ, PQPQP, PQRQP, PQSQP | 170 | (12) | (34) | (13) | (24) | Q, S, RPR, RSR, RQSP, RQPQR, RQRQR, RQSQSQR, RQSRSQR |
| 50 | (12) | (13) | (12) | (14)(23) | SPS, SRS, SQP, SQRQS, QSQS, RSQS, SQSPSQS | 171 | (12) | (13) | (13) | (14) | RS, RPR, RQR, PQP, PRP, PSP, QPQ, QRQ, QSQ |
| 51 | (12) | (13) | (14)(23) | (12) | PS, SQS, SRPRS, SRQRS, SRSRS, SPRPS, SPQPQPS, SPQRQPS, SPQSQPS | 172 | (12) | (13) | (14) | (13) | RQ, SPS, SQS, SRS, PQP, PRP, PSP, QPQ, QSQ |
| 52 | (12) | (14)(23) | (13) | (12) | SP, QPQ, QRQ, QSQ, PQR, PRP, RPR, RSR | 173 | (12) | (14) | (13) | (13) | RS, RPR, RQR, PQP, PRQ, PSP, QPQ, QSQ |
| 53 | (12) | (13) | (13) | (12)(34) | QR, RPR, RSPSR, RSQSR, RSRSR, RQSPQR, RQPQPQR, RQPRPQR | 174 | (23) | (24) | (14) | (1) | P, Q, S, RSQR, RPQPR, RPRPR, RPSPR, RQPQR, RQRQR |
| 54 | (12) | (13) | (12)(34) | (13) | QS, QPQ, QRPRQ, QRQRQ, QRSP, QRSR, QRSQSRQ | 175 | (23) | (12) | (1) | (24) | P, R, S, QPQ, QRQ, QSQPSQ, QSRSQ, QSPRPSQ, QSPSPSQ |
| 55 | (12) | (12)(34) | (13) | (13) | RS, SPS, SQPQS, SQRQS, SQSQS, SRQRS, SRPQPRS, SRPRPRS, SRPSPRS | 176 | (23) | (1) | (24) | (13) | P, Q, R, SQS, SRS, SPQPS, SPSPS, SPRPRPS, SPRQRPS, SPRSRPS |
| 56 | (12) | (13) | (13) | (13)(24) | RQ, SQ, PQP, PRP, QPQ, PSPSP, PSQSP, PSRSP | 177 | (23) | (12) | (24) | (1) | P, R, S, QPQ, QSQ, QRPRQ, QRSQRQ, QRQPQRQ, QRQRQRQ |
| 57 | (12) | (13) | (13)(24) | (13) | SQS, PS, SRPRS, SRQRS, SRSPSRS, SRSQPS, SRSRPS | 178 | (23) | (12) | (23) | (34) | P, R, S, QRPQ, QSQ, QPSQPQ, QPQPQPQ, QPQRQPQ |
| 58 | (12) | (13)(24) | (13) | (13) | RQ, SQ, PRP, PSP, QPQ, PQPQP, PQRQP, PQSQP | 179 | (23) | (12) | (34) | (24) | P, R, S, QRQ, QSQ, QPRPQ, QPSQPQ, QPQPQPQ, QPQRQPQ |
| 59 | (12) | (13) | (13) | (14)(23) | RQ, SPS, SQS, SRS, PQP, PRP, PSQ, QPQ | 180 | (23) | (12) | (24) | (34) | P, R, S, QPQ, QRQ, QSQSQ, QSRPSQ, QSPQPSQ, QSPSPSQ |
| 60 | (12) | (13) | (14)(23) | (13) | QS, SPQPQPS, SPQRPS, SPQSQPS, SPSPS, SRQS, SQPQS | 181 | (23) | (12) | (23) | (24) | P, R, S, QPQ, QRQ, QSQPSQ, QSRPSQ, QSPSPSQ |
| 61 | (12) | (14)(23) | (13) | (13) | SR, QPQ, QRQ, QSQ, PQR, PRP, PSP, RPR | 182 | (23) | (12) | (24) | (23) | P, R, S, QPQ, QSQ, QRQPRQ, QRSRQ, QRPRPRQ, QRPSPRQ |
| 62 | (12) | (13) | (14) | (12)(34) | RQP, PRP, PSP, PQSPSQP, PQSQSQP, PQSRSQP, PQPQP, PQRQRQP, PQRSRQP | 183 | (12)(34) | (12) | (13) | (1) | S, QP, RQR, RSR, PRP, PSP, RPQPR, RPRPR, RPSPR |
| 63 | (12) | (13) | (12)(34) | (14) | QP, RP, PSP, SPS, SQS, SRPRS, SRQRS, SRSRS | 184 | (12)(34) | (12) | (1) | (13) | R, RSR, RPQPR, RPRPR, RPSPR, RQRQR, RQSPQR, RQPQPQR, RQPRPQR |
| 64 | (12) | (12)(34) | (13) | (14) | QPQ, QRQ, QSR, QSQPQSQ, QSQRQSQ, QSQSQSQ, QSPSQ, QSRPRSQ, QSRSRSQ | 185 | (12)(34) | (1) | (12) | (13) | Q, RP, SQS, SRS, PQP, PSP, SPQPS, SPRPS, SPSPS |
| 65 | (12) | (13) | (14) | (13)(24) | SQ, RPR, RQR, PSR, QPQ, QRQ, RSQSR, RSRSR | 186 | (12)(34) | (13) | (12) | (1) | S, SRS, SPQPS, SPRPS, SPSPS, SQRPQS, SQSPQS, SQPQPQS |
| 66 | (12) | (13) | (13)(24) | (14) | RQ, QPQ, QSQ, SPS, SQS, SRP, SRQRS, SRSRS | 187 | (12)(34) | (13) | (1) | (12) | R, RSR, RPQPR, RPRPR, RPSPR, RQRPQR, RQSQR, RQPQPQR, RQPSPQR |
| 67 | (12) | (13)(24) | (13) | (14) | SPQ, QSQ, QPRPQ, QRSPQ, QPQPQPQ, QPQRQPQ, QPQSQPQ, QPSPSPQ | 188 | (12)(34) | (1) | (13) | (12) | Q, SP, RQR, RSR, PQP, PRP, RPQPR, RPRPR, RPSPR |
| 68 | (12) | (13) | (14) | (14)(23) | SQP, SPS, SQSPSQS, SQSQSQS, SQSRSQS, SQRQS, SQPQPQS, SQPRS | 189 | (12)(34) | (13) | (14) | (1) | S, PQP, PRP, PSP, QPR, QRQ, QSR, RQR |
| 69 | (12) | (13) | (14)(23) | (14) | QS, SRS, SPQPS, SPSPS, SQPQS, SQRQS, SPRPRPS, SPRQRPS, SPRSRPS | 190 | (12)(34) | (13) | (1) | (14) | R, PQP, PSP, PRP, QPS, QRQ, QSQ, SQS, SRS |
| 70 | (12) | (14)(23) | (13) | (14) | QRP, QPQ, QRPSQ, QRQPQRQ, QRQRQRQ, QRQSQRQ, QRSRQ, QRPRPRQ | 191 | (12)(34) | (1) | (13) | (14) | Q, PRP, PSP, PQP, RPS, RQR, RSR, SQS, SRS |
| 71 | (12) | (13) | (24) | (12)(34) | R, PS, PQP, PRSQ, PRPRP, PRQRP, PRSPSRP, PRSRSRP | 192 | (12)(34) | (12) | (12) | (13) | QP, QRQ, QSQ, QSQ, QPRPQ, QSQSQ, QSRSQ, QSPQPSQ, QSPRPSQ |
| 72 | (12) | (13) | (12)(34) | (24) | S, SPS, SQS, SRPRS, SRSQRS, SRQRPQRS, SRQPQPQRS, SRQPSPQRS | 193 | (12)(34) | (12) | (13) | (12) | QP, SP, PRP, PRP, PRQRP, PRSRP, PRPQPRP, PRPSPRP |
| 73 | (12) | (12)(34) | (13) | (24) | Q, S, RQP, PSQP, PQPQP, PRPRP, PRQRP, PRSRP | 194 | (12)(34) | (13) | (12) | (12) | SR, PQP, PRP, PSP, QPR, QRQ, QSQ, RQR |
| 74 | (12) | (13) | (24) | (13)(24) | R, RSR, RPRPRPR, RPRQRPR, RPRSPR, RPQPR, RQPQR, RQRQR, RQSQR | 195 | (12)(34) | (13) | (13) | (12) | PSP, PSPQPSP, PSPRPSP, PSPSPSP, PSQPRSP, PSQRQSP, PSQSQSP, PSRQRSP, PSRSRSP |
| 75 | (12) | (13) | (13)(24) | (24) | S, RPR, RQR, RSR, PQP, PRQ, PSP, QPQ, QSQ | 196 | (12)(34) | (13) | (12) | (13) | QS, RP, QRQ, PQP, PSP, QPQPQ, QPRPQ, QPSPQ |
| 76 | (12) | (13)(24) | (13) | (24) | S, RQ, PQSP, PRP, QPQ, QSQ, PQPQP, PQRQP | 197 | (12)(34) | (12) | (13) | (13) | RS, QP, RQR, PRP, PSP, RPQPR, RPRPR, RPSPR |
| 77 | (12) | (13) | (24) | (14)(23) | R, SPS, SQS, SRP, SRQRS, SSQSRS, SRSPSRS, SRSRSRS | 198 | (12)(34) | (12) | (13) | (14) | SRQ, QSQ, QRQRQ, QPSRQ, QRPQPRQ, QRPRPRQ, QRPSPRQ, QRSRSRQ |
| 78 | (12) | (13) | (14)(23) | (24) | S, SPS, SQS, SRQRS, SRPQSRS, SRPRPRS, SRPSPRS, SRSPSRS, SRSRSRS | 199 | (12)(34) | (14) | (12) | (13) | SQP, PSP, PQRQP, PRSQP, PQPQPQP, PQPRPQP, PQPSPQP, PQSQSQP |
| 79 | (12) | (14)(23) | (13) | (24) | S, SPS, SRQS, SQPQSQS, SQPRSQS, SQPSPQS, SQSPSQS | 200 | (12)(34) | (34) | (23) | (13) | Q, R, PQP, PRS, PSP, SQPS, SPRPS, SPSPS |
| 80 | (12) | (13) | (34) | (12)(34) | R, RSR, RPQPR, RPRPR, RPSPR, RQPQR, RQRPRQR, RQRQRQR, RQRSQR | 201 | (12)(34) | (12) | (13) | (24) | S, RQR, RSR, RPSP, RPSQ, RPQPR, RPRPR, RPSRSPR |
| 81 | (12) | (13) | (12)(34) | (34) | S, RP, QPQ, PQP, PSP, QRSQ, QRPRQ, QRQRQ | 202 | (12)(34) | (12) | (13) | (34) | SQS, SQPQPQS, SQPRPQS, SQPSQS, SQRQPRQS, SQRSPRQS, SQRPRPRQS |
| 82 | (12) | (12)(34) | (13) | (34) | S, SRS, SPSPS, SPQS, SPRPRPS, SPRQSRPS, SPRQPQRPS, SPRQRQRPS | 203 | (12)(34) | (12) | (34) | (13) | RQR, RQPQPQR, RQPRQR, RQPSPQR, RQSQPSQR, RQSRPSQR, RQSPSPSQR |
| 83 | (12) | (13) | (34) | (13)(24) | R, RQR,  RSPSR, RSQSR, RSRPRSR, RSRSQRSR, RPQRSR, RSRQRQRSR | 204 | (12)(34) | (13) | (12) | (24) | PRP, PRSRP, PRPQPRP, PRPRPRP, PRPSPQRP, PRPSQSPRP, PRPSRQRP, PRQSQRP |
| 84 | (12) | (13) | (13)(24) | (34) | S, SQS, SRPRS, SRQRS, SRSPSRS, SRSQPS, SRSRPS, SPSPS | 205 | (12)(34) | (13) | (24) | (12) | R, RSR, RPQPR, RPSPR, RQRQR, RPRPQR, RPRSQR, RPRQRPR |
| 85 | (12) | (13)(24) | (13) | (34) | QRQ, QSPR, RQSQ, QSPSQ, QSRSQ, QSQPQSQ, QSQSRQSQ | 206 | (12)(34) | (24) | (13) | (23) | Q, S, PRP, PSR, RPQP, RQR, PQRQP, PQSQP |
| 86 | (12) | (13) | (34) | (14)(23) | R, SPS, SQS, SRQ, SRPRS, PSRS, SRSQSRS, SRSRSRS | 207 | (12)(34) | (13) | (12) | (34) | RP, RQR, RPSR, RPQPSQPR, RPQRQPR, RPQPQPQPR, RPQPRPQPR |
| 87 | (12) | (13) | (14)(23) | (34) | S, RPR, RQR, QSR, PQP, PSP, RSPSR, RSRP | 208 | (12)(34) | (13) | (13) | (14) | RQ, PQP, PRP, PSP, QPS, QSQ, SQS, SRS |
| 88 | (12) | (14)(23) | (13) | (34) | S, QPQ, QRQ, RSQ, PRP, PSP, QSPSQ, QSQP | 209 | (12)(34) | (13) | (14) | (13) | SQ, PQP, PRP, PSP, QPR, QRQ, RQR, RSR |
| 89 | (12) | (13) | (24) | (1) | R, RSR, RQR, RPRPRPR, RPRQRPR, RPRSRPR, RPSQPR, RPQPQPR, RPQRQPR | 210 | (12)(34) | (14) | (13) | (13) | SR, PQP, PRP, QSP, RQR, PSPR, PSRSP |
| 90 | (12) | (13) | (1) | (24) | R, S, PQP, PRP, QPQ, QRQ, QSQ, PSPSP, PSQSP, PSRSP | 211 | (12)(34) | (13) | (14) | (12) | SP, SQS, SRS, SPQRPS, SPQSQPS, SPQPQPQPS, SPQPRPQPS, SPQPSPQPS |
| 91 | (12) | (1) | (13) | (24) | Q, S, PQP, PRP, RPR, RQR, RSR, PSPSP, PSQSP, PSRSP | 212 | (12)(34) | (23) | (24) | (12) | Q, R, SP, PQPRP, PQRQP, PQSQP, PRQRP, PRSRP |
| 92 | (12) | (13) | (34) | (1) | R, S, QPQ, QSQ, PQP, PRP, PSP, QRPRQ, QRQRQ, QRSRQ | 213 | (12)(34) | (23) | (24) | (13) | Q, R, PQS, PSP, SPRP, SRS, PRQRP, PRSRP |
| 93 | (12) | (13) | (1) | (34) | S, SQS, SRPRS, SRSRS, SRQPQRS, SRQSRQRS, SPRQRS, SRQRQRQRS | 214 | (12)(34) | (23) | (24) | (14) | Q, R, PRS, PSP, SPQP, SQS, PQRQP, PQSQP |
| 94 | (12) | (1) | (13) | (34) | Q, S, RPR, RQR, PQP, PRP, PSP, RSPSR, RSQSR, RSRSR | 215 | (12)(34) | (12) | (13) | (12)(34) | SP, SQS, SPQPS, SPRPS, SRSRS, SRQPRS, SRPRPRS, SRPSPRS |
| 95 | (12) | (34) | (13) | (1) | Q, S, RPR, RSR, PQP, PRP, PSP, RQPQR, RQRQR, RQSQR | 216 | (12)(34) | (12) | (12)(34) | (13) | RP, RQR, RPQPR, RSRSR, RPSPR, RSQPSR, RSPRPSR, RSPSPSR |
| 96 | (12) | (34) | (1) | (13) | Q, R, SPS, SRS, PQP, PRP, PSP, SQPQS, SQRQS, SQSQS | 217 | (12)(34) | (12)(34) | (12) | (13) | QP, RP, PSP, SRS, SPQS, SPRPS, SPSPS |
| 97 | (12) | (1) | (34) | (13) | Q, R, SPS, SQS, SRPRS, SRSRS, PQRS, SRQRQRS, SRQSQRS | 218 | (12)(34) | (12) | (13) | (13)(24) | RS, RQR, RPSP, RPSQ, RPQPR, RPRPR, RPSRSPR |
| 98 | (12) | (13) | (14) | (1) | S, SPS, SRS, SQSPSQS, SQSQSQS, SQSRSQS, SQRQS, SQPQPQS, SQPRPQS, SQPSPQS | 219 | (12)(34) | (12) | (13)(24) | (13) | RS, RQR, RPRP, RPRQ, RPQPR, RPSPR, RPRSRPR |
| 99 | (12) | (13) | (1) | (14) | QR, RPR, RSRPRSR, RSRQRSR, RSRSRSR, RSPSR, RSQSR, RQPQR, RQSQR | 220 | (12)(34) | (13)(24) | (12) | (13) | QS, QRQ, QPQP, QPQR, QPRPQ, QPSPQ, QPQSQPQ |
| 100 | (12) | (1) | (13) | (14) | Q, SPS, SQS, SRS, PQP, PRP, PSP, RPR, RQR, RSR | 221 | (12)(34) | (34) | (23) | (13)(24) | Q, R, PQP, PRS, SPSP, SQSP, PSRSP |
| 101 | (12) | (1) | (13) | (23) | Q, SPS, SQS, SRS, PQP, PRP, PSR, RPR, RQR | 222 | (12)(34) | (34) | (13)(24) | (23) | Q, S, PQP, PSR, RPRP, RQRP, PRSRP |
| 102 | (12) | (12)(34) | (12)(34) | (13) | QP, RP, PSP, SPS, SQPQS, SQRS, SQSQS | 223 | (12)(34) | (13)(24) | (34) | (23) | R, S, PRP, PSQ, QPQP, QRQP, PQSQP |
| 103 | (12) | (12)(34) | (13) | (12)(34) | QP, SP, PRP, RPR, RQSR, RQPQR, RQRQR | 224 | (12)(34) | (23) | (34) | (12)(34) | Q, R, SP, PRP, PQRPQP, PQSQP, PQPQPQP, PQPSPQP |
| 104 | (12) | (13) | (12)(34) | (12)(34) | RP, SP, PQP, QPQ, QRSQ, QRPRQ, QRQRQ | 225 | (12)(34) | (23) | (12)(34) | (34) | Q, RP, S, PSP, PQRPQP, PQSPQP, PQPQPQP |
| 105 | (12) | (34) | (13) | (12)(34) | Q, QSQ, QRPRQ, QPQPQ, QPRPQ, QPSPQ, QRQPQRQ, QRQRQRQ, QRQSRQ | 226 | (12)(34) | (12)(34) | (23) | (34) | QP, R, S, PSP, PRQPRP, PRSPRP, PRPRPRP |
| 106 | (12) | (34) | (12)(34) | (13) | Q, RPR, RQR, RSR, PQP, PRS, PSP, SPS, SQS | 227 | (12)(34) | (24) | (34) | (13)(24) | Q, R, PRP, PSQP, SPQP, SQS, SRQP |
| 107 | (12) | (12)(34) | (34) | (13) | QS, SPS, SRPRS, SRQRS, SRSRS, SQRQS, SQPQPQS, SQPRPQS, SQPSPQS | 228 | (12)(34) | (24) | (13)(24) | (34) | Q, S, PQR, PSP, RPRP, RSRP, PRQRP |
| 108 | (12) | (34) | (13) | (13)(24) | RPR, RSPSR, RSRSR, RSQRQSR, RSQSQSR, RQPQSR, PQSR, RSQPSPQSR | 229 | (12)(34) | (13)(24) | (24) | (34) | R, S, PRQP, PSP, QPQP, QRQ, QSQP |
| 109 | (12) | (34) | (13)(24) | (13) | Q, QSQ, QRPRQ, QRQRQ, QRSRQ, QPRPQ, QPSPQ, QPQPQPQ, QPQRQPQ, QPQSQPQ | 230 | (12)(34) | (24) | (34) | (14)(23) | Q, R, PQS, PRP, SPSP, SRSP, PSQSP |
| 110 | (12) | (13)(24) | (34) | (13) | R, RSR, RQPQR, RQSQR, RQRPRQR, RQRQPR, RQRSPR, RPRPR | 231 | (12)(34) | (24) | (14)(23) | (34) | Q, S, PRQP, PSP, RPQP, RQR, RSQP |
| 111 | (12) | (34) | (14) | (13)(24) | PS, PRP, PQPQP, PQRQP, PQSPSQP, PQSQPP, PQSR | 232 | (12)(34) | (14)(23) | (24) | (34) | R, S, PRQ, PSP, QPQP, QSQP, PQRQP |
| 112 | (12) | (34) | (13)(24) | (14) | Q, QPQ, QSQ, QRQPRQ, QRPRSRQ, QRPSPRQ, QRSPSRQ, QRSQSRQ | 233 | (12)(34) | (23) | (24) | (12)(34) | Q, R, SP, PQPRP, PQRQP, PQSRP, PRQRP |
| 113 | (12) | (13)(24) | (34) | (14) | R, QPQ, QSQ, SRQ, PRP, PSP, QRQP, QRPRQ | 234 | (12)(34) | (23) | (12)(34) | (24) | Q, RP, S, PQPSP, PQRSP, PQSQP, PSQSP |
| 114 | (12) | (23) | (12)(34) | (13) | P, Q, RQS, RSPR, SPS, SRS, RPQPR, RPRPR | 235 | (12)(34) | (12)(34) | (23) | (24) | QP, R, S, PRPSP, PRQSP, PRSRP, PSRSP |
| 115 | (12) | (13) | (24) | (12) | R, QPQ, QRQ, QSQ, PQS, PRP, PSP, SPS, SRS | 236 | (12)(34) | (24) | (23) | (13)(24) | Q, R, PRS, PSQP, SPQP, SQS, PQRQP |
| 116 | (12) | (13) | (12) | (24) | PR, RQR, RSPSR, RSQSR, RSRSR, RPSQPR, RPQPQPR, RPQRQPR | 237 | (12)(34) | (24) | (13)(24) | (23) | Q, S, PRQP, PSR, RPQP, RQR, PQSQP |
| 117 | (12) | (12) | (13) | (24) | *PQ, QRQ, QSPSQ, QSQSQ, QSRSQ, QPSPQ, QPRPRPQ, QPRQRPQ, QPRSRPQ* | 238 | (12)(34) | (13)(24) | (24) | (23) | *R, S, PRQP, PSQ, QPQP, QRQ, PQSQP* |
| 118 | (12) | (13) | (34) | (12) | *S, SRQS, SQPQS, SQSQS, SRPRS, SRQRS, SPSRS, SRSQSRS* | 239 | (12)(34) | (24) | (23) | (14)(23) | *Q, R, PQS, PRP, PSP, SRPS, SPQPS, SPSPS* |
| 119 | (12) | (13) | (12) | (34) | *S, RP, PQP, PSP, SQS, SQPQS, SQRQS, SQSPSQS, SQSRSQS* | 240 | (12)(34) | (24) | (14)(23) | (23) | *Q, S, PQR, PSRP, RPRP, RSR, PRQRP* |
| 120 | (12) | (12) | (13) | (34) | *S, QP, PRP, PSP, SRS, SRPRS, SRQRS, SRSPSRS, SRSQSRS* | 241 | (12)(34) | (14)(23) | (24) | (23) | *R, S, PRQ, PSQP, QPQP, QSQ, PQRQP* |
| 121 | (12) | (34) | (13) | (12) | *PS, SQS, SPQPS, SPRPS, SQRPRQS, SQRQS, SQRSRQS, SQRQPQRQS, SQRQSQRQS* | | | | | | |

Table 5. Colorings that will give rise to index 4 subgroups of $H$, where $\pi(H) \cong S_4$.

| index | *PQ* | *QR* | *RS* | *PR* | *PS* | *QS* | generators |
|---|---|---|---|---|---|---|---|
| 2 | (1) | (1) | (12) | (1) | (12) | (12) | *QP, RP, SRSP* |
| 3 | (123) | (132) | (23) | (1) | (23) | (13) | *RP, SP, QRPQ* |
| 4 | (234) | (143) | (13) | (123) | (23) | (34) | *QP, SP, RSRP* |

Table 6. Low index subgroups of a tetrahedron Kleinian group $K$ realized by the tetrahedron $t_{10}$.

| **Coxeter tetrahedron** | **Number of Subgroups of Tetrahedron Groups** | | | **Number of Subgroups of Tetrahedron Kleinian Groups** | | |
|---|---|---|---|---|---|---|
| | Index 2 | Index 3 | Index 4 | Index 2 | Index 3 | Index 4 |
| $t_1$ = [3, 5, 2, 3, 2, 2] | 3 | 1 | 1 | 1 | 1 | 0 |
| $t_2$ = [5, 3, 5, 2, 2, 2] | 1 | 0 | 0 | 0 | 0 | 0 |
| $t_3$ = [3, 5, 3, 2, 2, 3] | 1 | 0 | 0 | 0 | 0 | 0 |
| $t_4$ = [3, 5, 3, 2, 2, 2] | 1 | 0 | 0 | 0 | 0 | 0 |
| $t_5$ = [3, 4, 3, 2, 2, 3] | 1 | 0 | 0 | 0 | 0 | 0 |
| $t_6$ = [3, 5, 3, 2, 2, 4] | 1 | 0 | 0 | 0 | 0 | 0 |
| $t_7$ = [4, 3, 5, 2, 2, 2] | 3 | 0 | 1 | 1 | 0 | 0 |
| $t_8$ = [4, 3, 4, 2, 2, 3] | 3 | 1 | 4 | 1 | 1 | 3 |
| $t_9$ = [5, 3, 5, 2, 2, 3] | 1 | 0 | 0 | 0 | 0 | 0 |
| $t_{10}$ = [3, 3, 6, 2, 2, 2]; 1 | 3 | 1 | 2 | 1 | 1 | 1 |
| $t_{11}$ = [4, 4, 3, 2, 2, 2]; 1 | 7 | 1 | 8 | 2 | 1 | 4 |
| $t_{12}$ = [3, 3, 3, 3, 2, 2]; 1 | 1 | 1 | 5 | 0 | 1 | 1 |
| $t_{13}$ = [4, 3, 6, 2, 2, 2]; 1 | 7 | 2 | 9 | 3 | 2 | 4 |
| $t_{14}$ = [3, 4, 2, 4, 2, 2]; 1 | 7 | 1 | 10 | 3 | 1 | 6 |
| $t_{15}$ = [3, 6, 3, 2, 2, 2]; 2 | 3 | 2 | 5 | 1 | 3 | 0 |
| $t_{16}$ = [5, 3, 6, 2, 2, 2]; 1 | 3 | 0 | 1 | 1 | 0 | 0 |
| $t_{17}$ = [4, 3, 3, 3, 2, 2]; 1 | 1 | 1 | 3 | 0 | 1 | 4 |
| $t_{18}$ = [3, 6, 2, 3, 2, 2]; 2 | 3 | 4 | 3 | 1 | 4 | 2 |
| $t_{19}$ = [4, 4, 4, 2, 2, 2]; 3 | 15 | 0 | 35 | 7 | 0 | 29 |
| $t_{20}$ = [6, 3, 6, 2, 2, 2]; 2 | 7 | 3 | 49 | 3 | 3 | 3 |
| $t_{21}$ = [3, 4, 4, 2, 2, 3]; 1 | 3 | 1 | 4 | 1 | 1 | 5 |
| $t_{22}$ = [5, 3, 3, 3, 2, 2]; 1 | 1 | 1 | 1 | 0 | 1 | 1 |
| $t_{23}$ = [3, 3, 6, 2, 2, 3]; 2 | 1 | 1 | 3 | 0 | 1 | 4 |
| $t_{24}$ = [3, 3, 3, 3, 2, 3]; 2 | 1 | 4 | 4 | 0 | 4 | 5 |
| $t_{25}$ = [4, 4, 2, 4, 2, 2]; 2 | 15 | 0 | 35 | 7 | 0 | 31 |
| $t_{26}$ = [6, 3, 3, 3, 2, 2]; 3 | 1 | 1 | 51 | 0 | 1 | 4 |
| $t_{27}$ = [4, 3, 6, 2, 2, 3]; 2 | 3 | 4 | 5 | 1 | 4 | 5 |
| $t_{28}$ = [4, 4, 4, 2, 2, 3]; 2 | 7 | 1 | 14 | 3 | 1 | 2 |
| $t_{29}$ = [5, 3, 6, 2, 2, 3]; 2 | 1 | 0 | 0 | 0 | 0 | 0 |

| **Coxeter tetrahedron** | **Number of Subgroups of Tetrahedron Groups** | | | **Number of Subgroups of Tetrahedron Kleinian Groups** | | |
|---|---|---|---|---|---|---|
| | Index 2 | Index 3 | Index 4 | Index 2 | Index 3 | Index 4 |
| $\mathtt{t}_{30}$ = [6, 3, 6, 2, 2, 3]; 4 | 3 | 5 | 6 | 1 | 5 | 6 |
| $\mathtt{t}_{31}$ = [4, 4, 4, 2, 2, 4]; 4 | 15 | 0 | 35 | 7 | 0 | 49 |
| $\mathtt{t}_{32}$ = [3, 3, 3, 3, 3, 3]; 4 | 1 | 13 | 86 | 0 | 13 | 6 |

Table 7. The number of index 2, 3, 4 subgroups of tetrahedron and tetrahedron Kleinian groups generated by Coxeter tetrahedra of finite volume (up to conjugacy).